\documentclass[11pt,a4,fleqn]{article}
\usepackage{amsmath}
\usepackage{amssymb}
\usepackage{theorem}
\usepackage{xspace}
\def\pth#1{\left(#1\right)}
\def\acc#1{\left\{#1\right\}}
\def\cro#1{\left[#1\right]}

\def\uu{\textrm{\mathversion{bold}$\mathbf{1}$\mathversion{normal}}}
\def\oo{\textrm{\mathversion{bold}$\mathbf{0}$\mathversion{normal}}}
\def\eb{\textrm{\mathversion{bold}$\mathbf{\beta}$\mathversion{normal}}}  
\def\el{\textrm{\mathversion{bold}$\mathbf{\lambda}$\mathversion{normal}}}

\def\eE{I\!\!E}
\def\eP{I\!\!P}
\def\e1{1\!\!1}

\def\ef{\overset{.}{f}}

\def\XX{\textrm{\mathversion{bold}$\mathbf{X}$\mathversion{normal}}}
\def\UU{\textrm{\mathversion{bold}$\mathbf{u}$\mathversion{normal}}}
\def\pp{\textrm{\mathversion{bold}$\mathbf{p}$\mathversion{normal}}}

\theoremstyle{plain}

\newtheorem{theorem}{Theorem}[section]

\newtheorem{lemma}{Lemma}[section]
\newtheorem{remark}{Remark}
\newtheorem{proposition}{Proposition}[section]

\newcommand{\beqn}{\begin{eqnarray*}}
\newcommand{\eeqn}{\end{eqnarray*}}
\def\ef{\textrm{\mathversion{bold}$\mathbf{\phi}$\mathversion{normal}}}  
\def\ee1{\textrm{\mathversion{bold}$\mathbf{\varepsilon}$\mathversion{normal}}}

\def\eo{\textrm{\mathversion{bold}$\mathbf{\omega}$\mathversion{normal}}}

\newcommand{\N}{\mathbb{N}}
\newcommand{\R}{\mathbb{R}}

\def\argmin{\mathop{\mathrm{arg\,min}}} 
\begin{document}
\title {{\bf Quantile regression in high-dimension with breaking }}
\author{GABRIELA CIUPERCA  \footnote{{\it email: Gabriela.Ciuperca@univ-lyon1.fr}}}
\maketitle
\begin{center}
{\it Universit\'e de Lyon, Universit\'e Lyon 1, 
CNRS, UMR 5208, Institut Camille Jordan, 
Bat.  Braconnier, 43, blvd du 11 novembre 1918, 
F - 69622 Villeurbanne Cedex, France,}\\

\end{center}


\begin{abstract}
The paper considers a linear regression model in high-dimension for which the predictive variables can change the influence on the response variable at unknown times (called change-points). Moreover, the particular case of the heavy-tailed errors is considered. In this case, least square method with LASSO or adaptive LASSO penalty can not be used since the theoretical assumptions do not occur or the estimators are not robust. Then, the quantile model with SCAD penalty or median regression with LASSO-type penalty allows, in the same time, to estimate the parameters on every segment and eliminate the irrelevant variables. We show that, for the two penalized estimation methods, the oracle properties is  not affected by the change-point estimation. Convergence rates of the estimators for the change-points and for the regression parameters, by the two methods are found.  Monte-Carlo simulations illustrate  the performance of the methods.    
\end{abstract}


\section{Introduction}
A model which changes at some observations is called a change-point model. The location of these changes (called also change-points, breaks, changes) may be known or unknown. In this paper, we consider a model with multiple change-points at unknown locations. Moreover, as very often in practice, for example in genetics, the response variable is studied function of a very large number of regressors. However,  only a small number of  regressors is going to influence the response variable. 
In recent years, change-point models and high-dimension regression have received much attention in the literature, most often in the case of a model with zero mean errors and bounded variance. A $L_1$ or adaptive $L_1$ penalty in the context of least squares model can be considered. We obtain then the popular method introduced by Tibshirani(1996) and called LASSO (Least Absolute Shrinkage and Selection Operator) method.
On the other hand, it is well known that, the presence of outliers in model  may cause a large error in a least squares estimator. This can happen especially when the error distribution is not Gaussian and distribution tail is large enough. The outliers can also create problems in the detection of the jumps. An alternative method is then the quantile estimation.\\
To be more precise, if the errors $(\varepsilon_i)_{1 \leq i \leq n}$ of the regression model are such that $\eP[\varepsilon_i < 0]=\tau$, then the $\tau$th quantile regression is considered, i.e. the regression parameters are found by minimizing the function $\rho_\tau(\varepsilon)=\sum^n_{i=1} \varepsilon_i [ \tau \e1_{\varepsilon_i>0}-(1-\tau) \e1_{\varepsilon_i \leq 0} ]$. The choice of  $\tau=1/2$ yields the median regression and the  $L_1$-estimator, also known as least absolute deviation (LAD) estimator. \\
 Moreover, when the model has a very large regressor variable number, a penalty is necessary to estimate simultaneously the parameters on every segment and to eliminate the irrelevant regressors without crossing every time by a hypothesis test. The SCAD (Smoothly Clipped Absolute Deviation) and LASSO penalties have the advantage of selection and parameter estimation.  It was established that these two methods  have the oracle properties in a model without change-points: the zero components of the true parameters are estimated (shrunk) as 0 with probability tending to 1 (also called sparsity property) and the nonzero components have an optimal estimation rate (furthermore they are  asymptotically normal). See Wu and Liu (2009) for the SCAD method in a  $\tau$th quantile regression and Xu and Ying (2010) for the LASSO-type method in a median regression, both  models without change-points. Recall also for a median regression in high dimension the paper of Wang(2013), where a $L_1$ penalized least absolute deviation method is considered, when the overall variable number is larger than the observation number. \\
 In a multiple change-point model, the break estimation could affects the estimator properties. This is the main interest of this paper. The difficulty to study a change-point model results first from the dependence of the model of two parameter type: the regression  parameters and the change-points. \\
A change-point linear model in high-dimension was also considered by Ciuperca(2013) but under stronger assumptions that the errors have mean zero and  bounded variance. An adaptive LASSO estimator was studied. It was proved that it has the oracle properties on each estimated segment. However, when the model contains outliers, the adaptive LASSO estimator may not  be robust and moreover, the observation number should be greater than the parameter number to be estimated.    \\
In the present work we restrict our attention to the quantile regression in high-dimension with multiple change-points when the classical conditions on the errors do not occur. The change-points and the regression parameters on each segment are first estimated by the SCAD method. After, for a median regression ($\tau =1/2$), these parameters are estimated by the LASSO-type method. The asymptotic and oracle properties of these estimators are studied.     We also carry out simulations to investigate the properties of the two proposed estimators.\\
The paper is organized as follows. The model and assumptions are introduced in Section 2. In Section 3, the SCAD estimator in a change-point model is proposed and its asymptotic behavior is studied. Next, LASSO-type estimator is given in Section 4. For both methods, the oracle properties and convergence rate  of the estimators are obtained. Section 5 reports some simulations results which illustrate the methods interest. In Section 6 we give the proofs of Theorems. Finally, Section 7 contains some lemmas which are useful to prove the main results. 

\section{Model and general notations}
In this section we introduce the models without and with change-points, general assumptions, notations. Some general results used required for the two estimation methods are given.\\
We consider the linear  model  without  change-points
\begin{equation}
\label{eq1}
Y_i=\XX^t_i \ef +\varepsilon_i, \qquad i=1, \cdots, n
\end{equation}
where the response variable $Y_i$ is an univariate random variable,  $\XX_i \in \R^p$ is a $p$-vector of regressors (covariates) and the $\varepsilon_i$ is the error. The errors $(\varepsilon_i)_{1 \leq i \leq n}$ are independent identically distributed (i.i.d) random variables. The regression parameters are  $\ef \in  \Gamma \subset \R^p$, with $\Gamma$ a compact set and $\ef^0$ true value (unknown) of the parameter $\ef$.  Contrary to the classic suppositions for a regression  model, we do not impose the condition that the mean of errors $\varepsilon_i$ is zero or that their variance is  bounded.\\
All throughout the paper, vector and matrices are written in bold face. \\
With regard to the errors $\varepsilon_i$ and the design $\XX_i$, we make the following assumptions:\\
\textbf{(A1)} Let  $f$ be the density of  $\varepsilon_i$ and $F$ its distribution function. We suppose that $f(0)>0$, $F(0)=\tau$, $|f(y)-f(0)| \leq c |y|^{1/2}$, for all $y$ in a neighborhood of 0. The quantile $\tau$ is a real number in the interval  $(0,1) $.\\
\textbf{(A2)} $(\XX_i)_{1 \leq i \leq n}$ is a deterministic sequence, such that $n^{-1}\sum^n_{i=1} \XX_i \XX^t_i $ converges, as $n \rightarrow \infty$, to a non negative definite matrix;\\
\textbf{(A3)} $(\XX_i)_{1 \leq i \leq n}$ uniformly bounded.\\
These conditions are typical for a quantile regression (see e.g. Koenker, 2005). These assumptions are also classic  conditions  for a model  estimated by LAD method: the first  condition is found in Babu(1989) and the last two in  Bai(1998).\\

It is of interest to note that by assumption (A1) we have $\eP[\varepsilon_i < 0]=\tau$, but the expectation $\eE[\varepsilon_i]$ cannot exist. A regression model (\ref{eq1}) with the errors $(\varepsilon_i)$ satisfying the condition $\eP[\varepsilon_i <0]=\tau$ is called quantile regression.
In order to estimate the unknown regression parameter  $\ef$, we consider the  function
\begin{equation}
\label{eq2}
\rho_\tau(r)=  r [ \tau \e1_{r>0}-(1-\tau) \e1_{r \leq 0} ]
\end{equation}
and the corresponding estimator
\begin{equation}
\label{eq3}
\hat \ef^{(\tau)}_n=\argmin_{\ef \in \Gamma} \sum^n_{i=1} \rho_\tau(Y_i-\XX^t_i \ef).
\end{equation}

In order to study the quantile regression and the estimator (\ref{eq3}), let be the random  processes 
\begin{equation}
\label{eq4}
\begin{array}{l}
G^{(\tau)}_{i}(\ef;\ef^0)=\rho_\tau(\varepsilon_i-\XX^t_i(\ef-\ef^0))-\rho_\tau(\varepsilon_i), \qquad {\cal G}^{(\tau)}_{n}(\ef;\ef^0)=\sum^n_{i=1} G^{(\tau)}_{i}(\ef;\ef^0) ,\\
D_i=(1-\tau) \e1_{\varepsilon_i \leq 0} - \tau \e1_{\varepsilon_i >0} , \qquad W_n=\sum^n_{i=1} D_i \XX^t_i, \\
R^{(\tau)}_{i}(\ef;\ef^0)=G^{(\tau)}_{i}(\ef;\ef^0)-D_i \XX^t_i (\ef-\ef^0) .
\end{array}
\end{equation}
Obviously $\eE[D_i]=0.$ The relation between ${\cal G}^{(\tau)}_{n}$ and $R^{(\tau)}_{i}$ is
\begin{equation}
\label{eq5}
{\cal G}^{(\tau)}_{n}(\ef;\ef^0)-\eE[ {\cal G}^{(\tau)}_{n}(\ef;\ef^0)] =\sum^n_{i=1} [R^{(\tau)}_{i}(\ef;\ef^0) - \eE[R^{(\tau)}_{i}(\ef;\ef^0) ]] +W_n (\ef-\ef^0).
\end{equation}
For the parameter regression vector $\ef$, we shall use the notation  $\ef=(\phi_{,1}, \cdots, \phi_{,p})$.\\
Throughout the paper, $C$ denotes a positives generic constant not dependent on $n$ which may take different values in different formula or even in different parts of the same formula. For a vector $\mathbf{v}=(v_1, \cdots, v_p)$ let us denote $| \mathbf{v} | =(|v_1|, \cdots, |v_p|)$ and $\frac{1}{ \mathbf{v}}= (\frac{1}{v_1}, \cdots, \frac{1}{v_p})$. On the other hand, $\|\mathbf{v} \|_2$ is the Euclidean norm and $\|\mathbf{v} \|_1=\sum^p_{i=1}|v_i|$ is the $L_1$ norm. All vectors are column and $\textbf{v}^t$ denotes the transposed of $\textbf{v}$.\\

For coherence, we try to use the same notations as in the paper of Wu and Liu(2009).  
By elementary calculations, we obtain, with the probability 1, $| R^{(\tau)}_{i}(\ef;\ef^0)| \leq | \XX_i^t(\ef-\ef^0)| \e1_{|\varepsilon_i| \leq |\XX^t_i(\ef - \ef^0)| }$.
This inequality, the definition of $G^{(\tau)}_i$ and the assumption (A3) allow to obtain
${\cal G}^{(\tau)}_{n}(\ef;\ef^0) \leq C n \| \ef-\ef^0\|_2$.
Following result proves that for every parameter $\ef$ and for every quantile order $\tau$, the process  $G^{(\tau)}_{i}(\ef; \ef^0)$ has positive expectation, indifferently of the design $\XX_i$.

\begin{proposition}
\label{inegE0}
Under assumption (A1), we have, for all $\ef \in \Gamma$, $\eE[G^{(\tau)}_{i}(\ef; \ef^0) ] \geq 0$.
\end{proposition}

\begin{remark}
\label{Lemma4 Bai} 
In Bai(1998), the behavior in a neighborhood of $\ef^0$  of the process ${\cal G}^{(\tau)}_n -\eE[{\cal G}^{(\tau)}_n]$ is obtained in the particular case $\tau=1/2$. By a similar demonstration, we can prove that the result holds in general, for any $\tau \in (0,1)$: let be a positive sequence $(c_n)$ such that  $c_n \rightarrow 0$ and $n c^2_n/ \log n\rightarrow \infty$. Under the assumptions (A1)-(A3), there exists a constant $C >0$ such that $\forall \epsilon >0$, $\eP \cro{\sup_{\| \ef-\ef^0 \|_2 \leq c_n} \left|\frac{1}{n c^2_n}   \cro{ {\cal G}^{(\tau)}_{n}(\ef;\ef^0)-\eE[{\cal G}^{(\tau)}_{n}(\ef;\ef^0)]}\right| \geq \epsilon} \leq \exp(-\epsilon^2 n c_n^2 C)$.   
\end{remark}

The proof sketch  of this remark is given at the end of Section 6.\\
As a consequence of this Remark, by the Borel-Cantelli lemma, we have  for any $\epsilon >0$, 
\begin{equation}
\label{as}
\limsup_{n \rightarrow \infty} \left(\sup_{\|\ef - \ef^0\|_2 \leq c_n} \left|\frac{1}{n c^2_n}   \pth{ {\cal G}^{(\tau)}_{n}(\ef;\ef^0)-\eE[{\cal G}^{(\tau)}_{n}(\ef;\ef^0)]}\right|   \right)  \leq \epsilon, \qquad a.s.
\end{equation}

It is well known that the estimator (\ref{eq3}) has all nonzero components. For estimations and choosing the regressors  simultaneously, penalized methods can be used: SCAD or LASSO-type. These estimation methods all become  more interesting for a model with $K$ change-points
\begin{equation}
\label{eq15}
Y_i=\XX^t_i\ef_1 \e1_{1 \leq i < l_1}+ \XX^t_i \ef_2 \e1_{l_1 \leq i < l_2} + \cdots +\XX^t_i \ef_{K+1} \e1_{l_K \leq i \leq n}+\varepsilon_i, \qquad i=1, \cdots , n, 
\end{equation}
where $\e1_{(.)}$ denotes the indicator function. \\
The model parameters  are the regression parameters $(\ef_1, \cdots , \ef_{K+1})$ and the change-points $(l_1, \cdots , l_K)$. The true values (unknown) are $(\ef_1^0,\cdots, \ef^0_{K+1} )$, $(l^0_1, \cdots, l^0_K)$, respectively. The observations $l_{r-1}+1, \cdots , l_r$ between two consecutive change-points will be called the $r$th segment (interval, phase). \\
Concern the distance between two consecutive change-points, we impose the assumption\\
\textbf{(A4)} $l_{r+1} - l_r \geq n^{3/4}$, for all $r=0,1, \cdots , K$.\\

In order to study the properties of the penalized estimators in a model with breaking, we need  corresponding results obtained without change-points when $\tau=1/2$: by Ciuperca(2011b), Xu and Ying(2010) and for a some $\tau \in (0,1)$ by  Wu and Liu(2009).  \\
In the next section we investigate theoretical properties of the smoothly clipped absolute deviation (SCAD) method in a change-point model.

\section{SCAD estimator}
We begin this section by recalling the SCAD estimator for the quantile regression model (\ref{eq1}) without change-points, introduced by Fan and Li(2001) and developed later by Wu and Liu(2009)
\begin{equation}
\label{eq13}
\hat \ef^{(\tau,\lambda)}_{n} \equiv \argmin_\ef \pth{ \sum^n_{i=1} \cro{ \rho_\tau(Y_i-\XX^t_i \ef) +\sum^p_{j=1} p_\lambda(|\phi_{,j}|)}}.
\end{equation}
 The penalty $p_\lambda(\phi_{,j})$ is defined by its first derivative
\begin{equation}
\label{eq14}
p'_\lambda(|\phi_{,j}|) \equiv \lambda \acc{\e1_{|\phi_{,j}| \leq \lambda} +\frac{(a \lambda -|\phi_{,j}|)_{+}}{(a-1) \lambda} \e1_{|\phi_{,j}| > \lambda} }, \qquad \forall j=1, \cdots, p,
\end{equation}
with $\lambda >0$, $a>2$ deterministic tuning parameters.
For real $x$ we use the notation $sgn(x)$ for the sign function $sgn(x)=\frac{x}{|x|}$ when $x \neq 0$ and $sgn(0)=0$. We also denote $x_+=\max\{0,x \}$.\\
In order to study the  estimator $\hat \ef^{(\tau,\lambda)}_{n}$, introduce the function
\[
G^{(\tau,\lambda)}_{i}(\ef;\ef^0)\equiv G^{(\tau)}_{i}(\ef;\ef^0) +[\pp_\lambda(|\ef|)- \pp_\lambda(| \ef^0|)]^t \uu_p , \qquad i=1, \cdots, n,
\]
with $\uu_p\equiv (1, \cdots, 1)$ a $p \times 1$ vector and $\pp_\lambda(\ef) \equiv (p_\lambda (\phi_{,1}), \cdots p_\lambda(\phi_{,p})) $ also $p \times 1$ vector, with $\ef=(\phi_{,1}, \cdots, \phi_{,p})$. 

For this purpose, we first give the Karush-Kuhn-Tucker  (KKT) conditions for the quantile model (\ref{eq1}) without change-points.\\
For the estimator $\hat \ef^{(\tau , \lambda)}_n$ given by the relation (\ref{eq13}), let us consider the index set   of the variables selected by the SCAD method
\[
{\cal A}_n \equiv \acc{ j; \hat \phi^{(\tau,\lambda)}_{n,j} \neq 0},
\]
with $\hat \phi^{(\tau,\lambda)}_{n,j}$ the $j$th component of $\hat \ef^{(\tau,\lambda)}_{n}$.\\

\begin{proposition}
\label{Proposition KKT}
For the estimator (\ref{eq13}), the KKT conditions are
\[
\begin{array}{l}
  \textrm{for } j \in {\cal A}_n: \tau \sum^n_{i=1} X_{ij} - \sum^n_{i=1}X_{ij} \e1_{Y_i < \XX^t_i \hat \ef^{(\tau, \lambda)}_n}  =n \lambda  sgn(  \hat \phi^{(\tau, \lambda)}_{n,j}) \acc{\e1_{|  \hat \phi^{(\tau, \lambda)}_{n,j}| \leq \lambda } +\frac{(a \lambda - | \hat \phi^{(\tau, \lambda)}_{n,j}|)_+}{(a-1) \lambda }\e1_{| \hat \phi^{(\tau, \lambda)}_{n,j}| > \lambda}  },  \\
\textrm{for } j \not \in {\cal A}_n: \left| \tau \sum^n_{i=1} X_{ij} - \sum^n_{i=1}X_{ij} \e1_{Y_i < \XX^t_i \hat \ef^{(\tau, \lambda)}_n} \right| \leq n \lambda.   
\end{array}
\]
\end{proposition}


For the model (\ref{eq15}), in order to study the SCAD estimators of the regression parameters $(\ef_1, \cdots, \ef_{K+1})$, and of the change-points $(l_1, \cdots, l_K)$, let us consider the function
\begin{equation}
\label{eq16}
S(l_1,\cdots, l_K) \equiv \sum^{K+1}_{r=1} \inf_{(\ef_1, \cdots, \ef_{K+1}) \in \Gamma^{K+1}} \sum^{l_r}_{i=l_{r-1}+1} \cro{\rho_\tau(Y_i - \XX^t_i \ef_r)+\pp_{\lambda;(l_{r-1};l_r)} (|\ef_r |)\uu_p }.
\end{equation}
In each interval $(l_{r-1},l_r)$ another penalty $\pp_{\lambda;(l_{r-1};l_r)}$ can be considered, with $l_0=1$ and $l_{K+1}=n$.
For simplicity of notation, we denote the penalty of  (\ref{eq16}) by $\pp_{\lambda;(l_{r-1};l_r)}$ 
for some $(l_{r-1},l_r)$ and by $\pp_{\lambda;(l^0_{r-1};l^0_r)}$ for the true change-points, but it is understood that the 
series $\lambda$ are in fact $\lambda_{(l_{r-1};l_r})$, $\lambda_{(l^0_{r-1};l^0_r)}$, respectively. For the interval $(1,\cdots, n)$, the tuning parameter $\lambda_{(0,n)}$ is  $\lambda_n$. 

We define the SCAD change-point estimator by
\begin{equation}
\label{eq17}
(\hat l^{(\tau,\lambda)}_1, \cdots, \hat l^{(\tau,\lambda)}_K) \equiv \argmin_{(l_1, \cdots, l_K) \in \R^K} S(l_1, \cdots, l_K) ,
\end{equation}
with the function $S$ defined by (\ref{eq16}). 
Between two consecutive change-points $l_{r-1}$ and $l_r$, the SCAD estimator of the corresponding regression parameter $\ef_r$ is
\[
 \hat \ef^{(\tau,\lambda)}_{(l_{r-1};l_r)} \equiv \argmin_{\ef_r} \sum^{l_r}_{i=l_{r-1}+1} \cro{\rho_\tau(Y_i - \XX^t_i \ef_r)+\pp_{\lambda;(l_{r-1};l_r)} (|\ef_r |)\uu_p }=\argmin_{\ef_r}\sum^{l_r}_{i=l_{r-1}+1} G^{(\tau,\lambda)}_i (\ef_r; \ef^0_r) .
\]
Then, the SCAD  regression  parameter estimator for the $r$th segment is obtained by considering for the change-points their corresponding estimators: $\hat \ef^{(\tau,\lambda)}_{(\hat l^{(\tau; \lambda)}_{r-1};\hat l^{(\tau; \lambda)}_r)}$. 
The following two theorems state the asymptotic behaviors of the estimators (\ref{eq17}) and of $ \hat \ef^{(\tau,\lambda)}_{(\hat l^{(\tau; \lambda)}_{r-1};\hat l^{(\tau; \lambda)}_r)}$. 
The first result gives the convergence rate of the change-point estimator.
\begin{theorem}
\label{Theorem1 LASSO}
Under the assumptions (A1)-(A4), with the tuning parameter $(\lambda_{(l_{r-1},l_r)})_{1 \leq r \leq K+1}$   a sequence, depending on $n$, converging to zero, $(l_r-l_{r-1})^{1/2}\lambda_{(l_{r-1},l_r)} \rightarrow \infty$ and  for a  deterministic sequence $(c_n)$, such that   $c_n \rightarrow 0$, $n c^2_n/ \log n\rightarrow \infty$ and $\lambda_n c^{-2}_n  \rightarrow 0$, as $n \rightarrow \infty$, then  we have
$\hat l^{(\tau, \lambda)}_r-l^0_r =O_{\eP}(1)$, for every $r=1, \cdots, K$.
\end{theorem}

\begin{remark}
We have following relations between the sequences $(\lambda_n)$ and $(c_n)$,  $\lambda_n \ll c^2_n \ll c_n$. Example of sequence $(c_n)$ that satisfies the conditions in the Theorem \ref{Theorem1 LASSO}: $c^2_n=\lambda_n \log n$, for any sequence $(\lambda_n)$ converging to zero and $n^{1/2}\lambda_n  \rightarrow \infty$, as $n \rightarrow \infty$. An example of tuning parameter sequence $(\lambda_n)$ is the following $\lambda_n=n^{-2/5}$.
\end{remark}

By the Theorem 1 of Wu and Liu(2009), for tuning sequence  $\lambda_{(l_{r-1},l_r)} $ converging to zero as $n \rightarrow \infty$, we have that the convergence rate of the estimators of $\ef$ in each  segment is of order $(l^0_r-l^0_{r-1})^{-1/2}$. Hence, taking into account Theorem \ref{Theorem1 LASSO}, we deduce that $\| \hat \ef^{(\tau,\lambda)}_{(\hat l^{(\tau, \lambda)}_{r-1};\hat l^{(\tau, \lambda)}_r)}- \ef^0_r\|_2=O_{\eP} \pth{\hat l^{(\tau, \lambda)}_r - \hat l^{(\tau, \lambda)}_{r-1}}^{-1/2}$, for every $r=1, \cdots, K+1$, with $\hat l^{(\tau,\lambda)}_0=1$ and $\hat l^{(\tau,\lambda)}_{K+1}=n$.\\

We suppose that for each interval we have that the matrix $(l_r-l_{r-1})^{-1} \sum^{l_r}_{i=l_{r-1}+1} \XX_i \XX^t_i$ converges to  $ \textbf{C}_r$, as $n \rightarrow \infty$, with $\textbf{C}_r$ a non-negative definite matrix, which can be singular. Let us denote by $\textbf{C}^0_r$ the limiting matrix for the true change-points $l^0_r$, $r=1, \cdots, K$. We also denote by $C^0_{r,kj}$ the $(k,j)$th component of matrix $\textbf{C}^0_r$.\\
The following result proves that on every segment, the SCAD estimator for the regression parameters has the oracle properties: nonzero parameters estimator on each estimated segment is asymptotically normal and zero parameters are shrunk directly to 0 with a probability converging to 1. 
Let us underline that the limiting distribution not depend on the penalty $\pp_\lambda$, but only of the quantile order  $\tau$. For that purpose, for each two consecutive true change-points $l^0_{r-1}$, $l^0_r$ consider the set with the index of nonzero components of the true regression parameters 
\begin{equation}
\label{A0r}
{\cal A}_{(l^0_{r-1},l^0_r)} \equiv \acc{ j; \phi^0_{r,j} \neq 0}=^{noted}  {\cal A}^0_r
\end{equation}
 and with the index of the nonzero components of the SCAD regression parameter estimator  ${\cal A}_{n;(l^0_{r-1}, l^0_r)} \equiv \acc{j; \hat \phi^{(\tau,\lambda)}_{(l^0_{r-1},l^0_r),j} \neq 0 }$. Consider also the similar index set when the change-points are estimated ${\cal A}_{n;(\hat l^{(\tau, \lambda)}_{r-1}, \hat l^{(\tau, \lambda)}_r)} \equiv $ $\acc{j; \hat \phi^{(\tau, \lambda)}_{(\hat l^{(\tau, \lambda)}_{r-1},\hat l^{(\tau, \lambda)}_r),j} \neq 0}$. We denoted by $\ef_{ {\cal A}^0_r}$ the sub-vector of $\ef$ containing the corresponding components of $ {\cal A}^0_r$ and by $q_r\equiv Card \{{\cal A}^0_r \}$ the true number of nonzero components in the  $r$th segment.\\

\begin{theorem}
\label{Theorem6 LASSOadapt}
Under the assumptions (A1)-A4),  the tuning parameter sequence $(\lambda_{(l_{r-1},l_r)})_{1 \leq r \leq K+1}$ on each interval $(l_{r-1},l_r)$ as in Theorem \ref{Theorem1 LASSO}, then   we have\\
(i) $(\hat l^{(\tau,\lambda)}_r-\hat l^{(\tau,\lambda)}_{r-1})^{1/2} \pth{\hat \ef^{(\tau,\lambda)}_{(\hat l^{(\tau; \lambda)}_{r-1};\hat l^{(\tau; \lambda)}_r)}-\ef^0_r}_{{\cal A}^0_r }=(l^0_r-l^0_{r-1})^{1/2}\pth{\hat \ef^{(\tau,\lambda)}_{(\hat l^{(\tau; \lambda)}_{r-1};\hat l^{(\tau; \lambda)}_r)}-\ef^0_r}_{{\cal A}^0_r }(1+o_{\eP}(1))$ $ \overset{{\cal L}} {\underset{n \rightarrow \infty}{\longrightarrow}} {\cal N} \pth{\oo,\tau(1-\tau)/f^{2}(0)(\Omega^0_r)^{-1} }$, where  $\Omega^0_r \equiv (C^0_{r,kj})_{k,j \in {\cal A}_{(l^0_{r-1},l^0_r)}}$ is a $q_r \times q_r$ matrix.\\
(ii) 
$
\lim_{n \rightarrow \infty} \eP \cro{{\cal A}_{n;(l^0_{r-1}, l^0_r)}={\cal A}_{n;(\hat l^{(\tau, \lambda)}_{r-1}, \hat l^{(\tau, \lambda)}_r)}= {\cal A}^0_r}=1
$.
\end{theorem}

\section{LASSO-type estimator}
An important theoretical fact is that, as Zou(2006) showed recently,  the oracle properties do not hold for the LASSO estimator. We have just seen that considering the SCAD method, the obtained estimators have this property. But the last method is difficult to put into practice with regard to numerical algorithms. Thus, Xu and Ying(2010) proposed, for   model (\ref{eq1}), that the tuning parameter  $\lambda$ change from one component to the other of the  parameter $\ef$.\\
In this section the median model ( $\tau=1/2$) is studied. \\ 

Let us first consider, the model (\ref{eq1})  without change-points, mentioned in Section 2. The parameter  $\ef$ is estimate by
\begin{equation}
\label{eq24}
\hat \ef^L_n=\argmin_{\ef} \pth{ \sum^n_{i=1} |Y_i- \XX^t_i \ef |+\el^t_n |\ef|}.
\end{equation}
Compared with the SCAD  method seen in the previous section, now, the tuning parameter $\el_n=(\lambda_{n,1}, \cdots, \lambda_{n,p} )$ is a random $p$-vector with different components. 
The fact that  $\el_n$ has different components, makes possible that the estimator $\hat \ef^L_n$ have the oracle property, obviously, choosing the components of $\el_n$ in a judicious way. \\
For the regression  model (\ref{eq1}) without change-points, and for the estimator  $\hat \ef^L_n$ given by (\ref{eq24}), consider the index set of estimator nonzero components  ${\cal A}^L_n \equiv \{ j; \hat \phi^L_{n,j} \neq 0 \}$ where $\hat \phi^L_{n,j}$ the $j$th component of $\hat \ef^L_n$. Similar to the Proposition \ref{Proposition KKT}, we obtain that the KKT  relations are in this case:\\
 $ - \sum^n_{i=1} X_{ij} \cdot sgn(Y_i-\XX^t_i \hat \ef^L_n)+\lambda_{n,j}\cdot sgn( \hat \phi^L_{n,j})=0$,   for all $j \in {\cal A}^L_n$,\\
 $ \left| \sum^n_{i=1} X_{ij}\cdot sgn(Y_i-\XX^t_i \hat \ef^L_n) \right| \leq \lambda_{n,j}$, for all $j \not \in {\cal A}^L_n$,\\ 
with  $\lambda_{n,j}$ the $j$th  component of  $\el_n$ and $ \hat \phi^L_{n,j}$ of $\hat \ef^L_n$. These results will be useful to prove the oracle properties for the LASSO-type estimators of the regression parameters on each segment, in a model with change-points.  \\

Consider now the change-point problem (\ref{eq15}), with $K$ (known) changes. For this estimation method,
the change-point estimator is
\[
(\hat l^L_1, \cdots, \hat l^L_K) \equiv \argmin_{(l_1, \cdots, l_K) \in \R^K} \sum^{K+1}_{r=1} \inf_{(\ef_1, \cdots, \ef_{K+1})} \sum^{l_r}_{i=l_{r-1}+1}\cro{|Y_i-\XX^t_i \ef_r| +\frac{\el^t_{n,(l_{r-1},l_r)}}{l_r-l_{r-1}} |\ef_r|}.
\]
The LASSO-type estimator of the regression parameters for the $r$th segment is $\hat \ef^L_{(\hat l^L_{r-1};\hat l^L_r)}$, for each $r=1, \cdots , K+1$, with $\hat l^L_0=1$ and $\hat l^L_{K+1}=n$.
Taking into account that a particular case ($\tau=1/2$) to the  quantile regression is considered, following processes are introduced
\begin{equation}
\label{eq25}
\begin{array}{l}
G^{(1/2)}_{i}(\ef,\ef^0) \equiv |\varepsilon_i-\XX^t_i(\ef-\ef^0) | - | \varepsilon_i|, \qquad i=1, \cdots, n \\
\eta^L_{i;(j_1,j_2)}(\ef, \ef^0) \equiv G^{(1/2)}_{i}(\ef,\ef^0)+ (j_2-j_1)^{-1} \el^t_{n;(j_1,j_2)}(|\ef| -|\ef^0 | ), \qquad i=j_1+1, \cdots, j_2
\end{array}
\end{equation}
with $0 \leq j_1 < j_2 \leq n$ and $\ef^0$  the true parameter. In the particular case $j_1=0$ and $j_2=n$, let us denote $\lambda_{n;(0,n)}$ by $\lambda_n$.\\
Observe that, since we will study the model (\ref{eq15}) with change-points, by the least absolute deviation method ($\tau=1/2$) with  LASSO-type penalty,  the related  results obtained when there is no penalty by Bai(1998), Ciuperca(2011b) are needed.\\

Following result yields that, even if the penalty is different, this estimator has the same convergence rate as the estimator obtained by the SCAD method. 
\begin{theorem}
\label{Theorem1 LASSObis} 
If the tuning parameter $\el_{n,(l_{r-1},l_r)}$ satisfies the conditions $\| \el^t_{n,(l_{r-1},l_r)} \|_2 \rightarrow \infty$, $(l_r-l_{r-1})^{-1/2} \| \el^t_{n,(l_{r-1},l_r)} \|_2 \overset{\eP} {\underset{n \rightarrow \infty}{\longrightarrow}} M \geq 0  $, under the assumptions (A1)-(A4), we have
$\hat l^L_r - l^0_r=O_{\eP}(1)$, for every $r=1, \cdots, K$.
\end{theorem}

Combining the Theorem \ref{Theorem1 LASSObis} and the $\sqrt{n}$-consistency of the parameter estimator in a model without change-points (see Theorem 2 of Xu and Ying, 2010), we have that the convergence rate of the  regression parameter  LASSO-type estimator on each segment is $\|\hat \ef^L_{(\hat l^L_{r-1},\hat l^L_r)} - \ef^0_r \|_2=(l^0_r-l^0_{r-1})^{-1/2}O_{\eP}(1)$, for  $r=1, \cdots , K+1$, with $l^0_0=1$ and $l^0_{K+1}=n$.\\

For this type of  method, the most important is to verify that if the oracle properties are preserved in a change-point model. The sparsity property is the most interesting and it risk to be influenced by the change-point estimation. We would like to point out that, due to a penalty different, the proof of this result differs from that for the SCAD estimator.

\begin{theorem}
\label{Theorem6 LASSOadaptbis}
Under the assumptions (A1)-(A4), with the tuning sequence $(\el_{n,(l_{r-1},l_r)})$ as in Theorem \ref{Theorem1 LASSObis} and the index set ${\cal A}^0_r $ defined by (\ref{A0r}), we have:\\
(i) $(\hat l^L_r-\hat l^L_{r-1})^{1/2} \pth{\hat \ef^L_{(\hat l^L_{r-1};\hat l^L_r)}-\ef^0_r}_{ {\cal A}^0_r }=(l^0_r-l^0_{r-1})^{1/2}\pth{\hat \ef^L_{(\hat l^L_{r-1};\hat l^L_r)}-\ef^0_r}_{ {\cal A}^0_r }(1+o_{\eP}(1))$ converges in distribution to the $p$-dimensional Gaussian vector $ {\cal N} \pth{\oo, 1/(4f^2(0))(\Omega^0_r)^{-1} }$, as $n \rightarrow \infty$.\\
(ii) 
$
\lim_{n \rightarrow \infty} \eP \cro{{\cal A}^L_{n;(l^0_{r-1}, l^0_r)}={\cal A}^L_{n;(\hat l^{L}_{r-1}, \hat l^{L}_r)}= {\cal A}^0_r}=1
$, 
with:  ${\cal A}^L_{n;(l^0_{r-1}, l^0_r)} \equiv \acc{j; \hat \phi^L_{(l^0_{r-1},l^0_r),j} \neq 0 }$ and ${\cal A}^L_{n;(\hat l^L_{r-1}, \hat l^L_r)} \equiv \acc{j; \hat \phi^L_{(\hat l^L_{r-1},\hat l^L_r),j} \neq 0}$.
\end{theorem}

It is worthwhile to mention that, if the same model (\ref{eq15}) is estimated by least squares, under certain conditions on design $(\XX_i)$, with a LASSO penalty, the sparsity property (i.e. the claim (ii) of the Theorem \ref{Theorem6 LASSOadaptbis}),  is not satisfied (see Zou, 2006).  Moreover, as the model contains change-points, this condition is more difficult to check on each interval that has random bounds. Then, an adaptive LASSO method can be considered downside to remedy this. But, it is necessary that in each segment $(l_{r-1}, l_r)$ parameter number is smaller than observation number $ l_{r-1}-l_r$. On the other hand, the adaptive LASSO for least squares method holds only under the assumptions that the errors have mean zero and bounded variance.\\

\textbf{An example of tuning random sequence $\el_{n,(l_{r-1},l_r)}$}: in each segment $(l_{r-1},l_r)$ the LAD estimator $\hat \ef^{(1/2)}_{(l_{r-1},l_r)}$ of $\ef_r$ is calculated by a corresponding relation to (\ref{eq3}), for $\tau=1/2$. Obtained estimators have all nonzero components and they have a convergence rate $v_{n,(l_{r-1},l_r)}$ to the true parameter, with $(l_r-l_{r-1}) v_{n,(l_{r-1},l_r)} \rightarrow \infty  $ (see Theorem 1 of Ciuperca, 2011b). Consider then  $\el_{n,(l_{r-1},l_r)}=\pth{ \frac{1}{|\hat \phi^{(1/2)}_{(l_{r-1},l_r),1}|}, \cdots , \frac{1}{|\hat \phi^{(1/2)}_{(l_{r-1},l_r),p}|}}$. 

\section{Simulation study}
We now give some simulation results. 
All simulations were performed using the R language. To calculate Least squares estimation the function \textit{lm} was used. While, for the quantile estimations, SCAD and LASSO-type, the function \textit{rq} of the package \textit{quantreq} were called. To compare these estimates when the classical conditions on the error distribution do not occur, we consider also  the adaptive LASSO estimation using  the function \textit{lqa} of the package \textit{lqa} and quantile estimation with LASSO penalty.\\
The number of phases is assumed to be known: the models contain two change-points (three phases).  We consider 10 latent variables $X_1, \cdots, X_{10}$ with $X_3 \sim {\cal N}(2,1)$, $X_4 \sim {\cal N}(4,1)$, $X_5 \sim {\cal N}(1,1)$ and $X_j \sim {\cal N}(0,1)$ for $j \in \{1,2,6,7,8,9,10 \}$.   The true values of the regression parameters (coefficients) on the three segments are respectively: $(1,0,4,0,-3,5,6,0,-1,0)$, $(0,3,-4,-3,0,1,2,-3,0,10)$, $(1,3,4,0,0,1,0,0,0,1)$. Three error patterns were considered: exponential, Cauchy and standard normal distributions. For the exponential errors, we generate a $n$-sample of distribution ${\cal E}xp(-1.5,1)$, with the density $\exp (-(x+1.5)) \e1_{x > -1.5}$. 
 For each model, we generated 500 Monte-Carlo random samples of size $n$, with $n=60$ or $n=200$. The percentage of zero coefficients correctly estimated to zero(true 0) and the percentage of nonzero coefficients estimated to zero(false 0) are computed (see Tables \ref{Tabl1}-\ref{Tabl6}) by least squares(LS), quantile(QUANT), quantile with LASSO penalty(QLASSO), SCAD,  LASSO-type and adaptive LASSO  methods. The reader can find in the paper of Ciuperca(2013) more details on the adaptive LASSO method in a change-point model.  The adaptive LASSO estimators of the change-points and of the regression parameters are the minimizers of the following penalized sum $\sum^{K+1}_{r=1} [\sum^{l_r}_{i=l_{r-1}+1} (Y_i+\XX^t_i \ef)^2 + \lambda_{n;(l_{r-1},l_r)}  \hat \eo_{(l_{r-1},l_r)} |\ef|]$, where the adaptive penalty $p$-vector $\hat \eo_{(l_{r-1},l_r)}$ is considered here that $| \hat \ef^{LS}_{(l_{r-1},l_r)}|^{-9/40}$. Let us specify that $\hat \ef^{LS}_{(l_{r-1},l_r)}$ is the LS estimator of $\ef$ calculated between $l_{r-1}$ and $l_r$. Recall also that the adaptive LASSO estimator of the regression parameters has the oracle properties under the assumptions for the errors $\varepsilon$ that $\eE[\varepsilon]=0$ and $\eE[\varepsilon^2] < \infty$.  For the quantile method with the LASSO  penalty, the sum $\sum_{i=l_{r-1}+1}^{l_r} \rho_\tau(Y_i+\XX^t_i \ef)$ is penalized with $\lambda_{n;(l_{r-1},l_r)} |\ef|$. 
 The tuning parameters $\lambda_{n;(l_{r-1},l_r)}$ are  $\log(l_r- l_{r-1}) \cdot \uu_p$ for the quantile estimation with LASSO penalty, $(l_r-l_{r-1})^{-2/5}$ for SCAD and $(l_r-l_{r-1})^{2/5}$ for adaptive LASSO methods.  For the  LASSO-type method, the tuning parameter is $(l_r-l_{r-1})^{2/5}\cdot 1/\hat \ef^{QLASSO}_{(l_{r-1},l_r)}$, where $\hat \ef^{QLASSO}_{(l_{r-1},l_r)}$ is the corresponding estimate by the quantile method (for  the index quantile of the errors equal to $\tau$) with LASSO penalty. 
  Since the asymptotic distribution of the change-points estimators can not be symmetric,  in each table we also give the median of the change-point estimations. Because the results by the SCAD method are  poorer than  by the LASSO-type method, and also because there may be convergence problem (the function \textit{rq} not responding), in Tables \ref{Tabl4}-\ref{Tabl6} the SCAD estimator is not considered.\\
The outliers of the errors do not affect the precision of the change-point estimations, by all six methods, while the sparsity property of the QLASSO and adaptive LASSO are affected. More specifically, when $n$ is large enough ($n=200$) and the errors $\varepsilon$ are normal, then $\eE[\varepsilon]=0$ and $\eE[\varepsilon^2] < \infty $, the two methods, adaptive LASSO and LASSO-type, give the same (very satisfactory) sparsity results (as Ciuperca, 2013, also indicates, for the adaptive LASSO estimators in a change-point model). When $n$ or number of observations in a segment is small, the LASSO-type method is better than the adaptive LASSO method, in terms of detection of irrelevant regressors (true and false zeros). If the errors are ${\cal E}(-1.5,1)$, then $\eE[\varepsilon] \neq 0$, the results for LASSO-type are relatively better than for adaptive LASSO method (see Tables \ref{Tabl1} and \ref{Tabl4}). This difference is accentuated when the moments of errors don't exist, $\varepsilon \sim Cauchy$ (see Tables \ref{Tabl3} and \ref{Tabl6}). Since LASSO-type and adaptive LASSO methods gave the best results, we calculate the average of estimation error $\| \hat \ef -\ef^0\|_1$ in each segment, for the index corresponding to the true values different to  zero,  over 500 simulations for different error distributions, for $n=200$ and two change-points  $l^0_1=30$, $l^0_2=100$ (see Table \ref{Tabl7}). For Gaussian and exponential distributions these two estimation methods, yield similar results. On the other hand, for Cauchy distribution, the obtained estimations by adaptive LASSO method are biased.  \\
  In conclusion, the LASSO-type method provides very satisfactory estimations in any case even for small sample size. The only less favorable result is obtained for $n=60$ when the errors are exponential. The percentage of false zero is large enough.

\section{Proofs of Theorems and Propositions}
In order to simplify the proofs of theorems and propositions, we give in this section their demonstrations and in Section 6 some lemmas and their proofs which will useful.\\

\noindent {\bf Proof of Proposition \ref{inegE0}}\\
Let us consider the notations $h_i(\ef)=\XX^t_i(\ef^0-\ef)$ and $F(x)$ for the distribution function of  $\varepsilon_i$. By definition 
$
\eE[G^{(\tau)}_{i}(\ef;\ef^0)]= \int_{\R} [\rho_\tau(x+h_i(\ef))-\rho_\tau(x) ] dF(x).
$
 Using $F(0)=\tau$, a simple algebraic computation gives 
\begin{equation}
 \label{eq7}
 \eE[G^{(\tau)}_{i}(\ef;\ef^0) ]=\int^{-h_i(\ef)}_0 [|h_i(\ef)|-x] dF(x), \qquad \textrm{if }h_i(\ef) <0.
\end{equation}
and
\begin{equation}
\label{eq8}
 \eE[G^{(\tau)}_{i}(\ef;\ef^0) ]=\int_{-h_i(\ef)}^0 [|h_i(\ef)|+x] dF(x), \qquad \textrm{if }h_i(\ef) \geq 0.
\end{equation}
Taking into account the  relations (\ref{eq7}) and (\ref{eq8}) we can write
\[ 
\eE[G^{(\tau)}_{i}(\ef;\ef^0)]
\geq \e1_{h_i(\ef) \geq 0}\int_{-\frac{h_i(\ef)}{2}}^0 [|h_i(\ef)|+x] dF(x) +\e1_{h_i(\ef) <0} \int^{-\frac{h_i(\ef)}{2}}_0 [|h_i(\ef)|-x] dF(x)
\]
\[
\geq \e1_{h_i(\ef) \geq 0} \frac{|h_i(\ef)|}{2} \int_{-\frac{h_i(\ef)}{2}}^0 dF(x) +\e1_{h_i(\ef) <0} \frac{|h_i(\ef)|}{2}  \int^{-\frac{h_i(\ef)}{2}}_0 dF(x)
\]
\[
=\frac{|h_i(\ef)|}{2} \cro{\e1_{h_i(\ef) \geq 0} [F(0)-F(-\frac{h_i(\ef)}{2})]+\e1_{h_i(\ef) <0}  [F(-\frac{h_i(\ef)}{2})-F(0)]} \geq 0.
\]
Hence, $\eE[G^{(\tau)}_{i}(\ef;\ef^0) ] \geq 0$ , for all $i=1, \cdots , n$,  $\ef \in \Gamma $.
\hspace*{\fill}$\Diamond$ \\

\subsection{For SCAD estimator}
\noindent {\bf Proof of Proposition \ref{Proposition KKT}}\\
\underline{ If $j \in {\cal A}_n$.} According to the definition (\ref{eq13}), the SCAD estimator of $\ef$ is the solution of the following equation
$
0= \sum^n_{i=1}\frac{\partial G^{(\tau,\lambda)}_{i}(\hat \ef^{(\tau,\lambda)}_{n})}{\partial \phi_{,j}} 
$
\[= \sum^n_{i=1}\pth{- \tau X_{ij} +X_{ij} \e1_{Y_i < \XX^t_i \hat \ef^{(\tau, \lambda)}_n} + \lambda \acc{\e1_{|  \hat \phi^{(\tau, \lambda)}_{n,j}| \leq \lambda } \cdot sgn(  \hat \phi^{(\tau, \lambda)}_{n,j}) +\frac{(a \lambda - | \hat \phi^{(\tau, \lambda)}_{n,j}|)_+}{(a-1) \lambda } sgn( \hat \phi^{(\tau, \lambda)}_{n,j}) \e1_{| \hat \phi^{(\tau, \lambda)}_{n,j}| > \lambda} }}.
\]
We obtain
\[
 \tau \sum^n_{i=1} X_{ij} - \sum^n_{i=1}X_{ij} \e1_{Y_i < \XX^t_i \hat \ef^{(\tau, \lambda)}_n} = n \lambda \cdot sgn(  \hat \phi^{(\tau, \lambda)}_{n,j}) \acc{\e1_{|  \hat \phi^{(\tau, \lambda)}_{n,j}| \leq \lambda } +\frac{(a \lambda - | \hat \phi^{(\tau, \lambda)}_{n,j}|)_+}{(a-1) \lambda }\e1_{| \hat \phi^{(\tau, \lambda)}_{n,j}| > \lambda}  }.
\]
\underline{ If $j \not  \in  {\cal A}_n$.}  In this case
$
0 \in \sum^n_{i=1} \frac{\partial G^{(\tau, \lambda)}_{i,n} (\hat \phi^{(\tau, \lambda)}_{n,j})}{\partial \phi_{,j}}=- \tau \sum^n_{i=1} X_{ij} +\sum^n_{i=1}X_{ij} \e1_{Y_i < \XX^t_i \hat \ef^{(\tau, \lambda)}_n} + \sum^n_{i=1}  p'_\lambda(|\hat \phi^{(\tau, \lambda)}_{n,j}|)$. Since 
$
p'_\lambda(|\phi_{,j}|) \equiv \frac{\partial p_\lambda(|\phi_{,j}|)}{\partial \phi_{,j}}=\lambda \cdot sgn(\phi_{,j}) \e1_{|\phi_{,j}| \leq \lambda}+\frac{(a \lambda -|\phi_{,j}|)_+ \cdot sgn(\phi_{,j})}{a - 1} \e1_{|\phi_{,j}| > \lambda},$
it follows that 
$
0 \in \sum^n_{i=1} \frac{\partial G^{(\tau, \lambda)}_{i,n} (\hat \phi^{(\tau, \lambda)}_{n,j})}{\partial \phi_{,j}}=- \tau \sum^n_{i=1} X_{ij} +\sum^n_{i=1}X_{ij} \e1_{Y_i < \XX^t_i \hat \ef^{(\tau, \lambda)}_n} + n \lambda \cdot [-1,1]$. 
Then, $| \tau \sum^n_{i=1} X_{ij} - \sum^n_{i=1}X_{ij} \e1_{Y_i < \XX^t_i \hat \ef^{(\tau, \lambda)}_n} | \leq n \lambda $. 
\hspace*{\fill}$\Diamond$ \\

\noindent {\bf Proof of Theorem \ref{Theorem1 LASSO}}\\
The proof is similar to that of Theorem \ref{Theorem1 LASSObis}. It is omitted. The Lemmas \ref{Lemma1 LASSO} and \ref{Lemma3 LASSO} stated in Section 6 are needed.
\hspace*{\fill}$\Diamond$ \\

\noindent {\bf Proof of Theorem \ref{Theorem6 LASSOadapt}}\\
(i) The statement results from  Theorem \ref{Theorem1 LASSO} together with  Theorem 2(b) of Wu and Liu (2009).\\
(ii) By the  Theorem 2(a) of Wu and Liu(2009), we have: $\lim_{n \rightarrow \infty} \eP \cro{{\cal A}_{n;(l^0_{r-1}, l^0_r)}={\cal A}_{(l^0_{r-1},l^0_r)}}=1$. The asymptotic normality of the estimators implies: for all $k \in {\cal A}_{(l^0_{r-1},l^0_r)} $, we have $\phi^0_{r,k} - \hat \phi^{(\tau, \lambda)}_{(\hat l^{(\tau, \lambda)}_{r-1}, \hat l^{(\tau, \lambda)}_r),k} \overset{{\eP}} {\underset{n \rightarrow \infty}{\longrightarrow}} 0$. It follows that
\begin{equation}
\label{eq21}
\lim_{n \rightarrow \infty} \eP \cro{ {\cal A}_{n;(\hat l^{(\tau, \lambda)}_{r-1}, \hat l^{(\tau, \lambda)}_r)} \supseteq  {\cal A}_{(l^0_{r-1},l^0_r)}}=1.
\end{equation}
By similar arguments that for the  Lemma 1 of Wu and Liu(2009), we prove that
\begin{equation}
\label{eq22}
\eP \cro{\exists k \in \{1, \cdots, p \}, k \not \in {\cal A}_{(l^0_{r-1},l^0_r)}, k \in {\cal A}_{n;(\hat l^{(\tau, \lambda)}_{r-1}, \hat l^{(\tau, \lambda)}_r)} } {\underset{n \rightarrow \infty}{\longrightarrow}} 0.
\end{equation}
The Theorem results from the relations (\ref{eq21}) and (\ref{eq22}).
\hspace*{\fill}$\Diamond$ \\

\subsection{For LASSO-type estimator}

\noindent {\bf Proof of Theorem \ref{Theorem1 LASSObis}}\\
The proof has three steps. First, we show that the all SCAD estimators of the change-points are to a smaller distance than $n^{1/2}$ from the corresponding true value. Then, for each true change-point $l^0_r$, with $r \in \{1, \cdots, K \}$, we consider the function $S$ given by (\ref{eq16}), but calculated on the change-points $l_1, \cdots, l_K, l^0_1, \cdots, l^0_{r-1}, l^0_r-[n^\alpha], l^0_r+[n^\alpha], l^0_{r+1}, \cdots , l^0_K$, with $\alpha \in (1/2, 1)$. For the penalized sums involving observations between $l^0_{t-1}$ and $l^0_t$, for $t \in \{1, \cdots, r-1, r+1, \cdots , K \}$, consider the change-points $k_{1,t} < \cdots < k_{J(t),t} \equiv \{l_1, \cdots, l_K \} \cap \{j; l^0_{r-1} < j \leq l^0_r \}$. Then for each $t \in \{1, \cdots, r-1, r+1, \cdots , K \}$, we have 
\[
0 \geq \sum^{J(t)+1}_{j=1} \min_{\ef_j \in \Gamma} \cro{\sum^{k_{j,t}}_{i=k_{j-1},t+1} |\varepsilon_i-\XX^t_i(\ef_j-\ef^0_t) |+ \el^t_{n;(k_{j-1,t};k_{j,t})} \cdot |\ef_j|}
\]
\[ - \sum^{J(t)+1}_{j=1} \cro{ \sum^{k_{j,t}}_{i=k_{j-1},t+1}| \varepsilon_i | + \el^t_{n;(k_{j-1,t};k_{j,t})} |\ef^0_j|} \geq -2(K+1) \sup_{1 \leq l < j \leq n} \left| \inf_{\ef} \cro{ \sum^j_{i=l+1} \eta^L_{i;(l,j)}(\ef,\ef^0) } \right|, 
\]
which is, using the Lemma \ref{Lemma1 LASSObis}, $-O_{\eP}(\max(n^\alpha,\el_n))$, with $\alpha \in (1/2, 1)$.    The rest of proof  is similar to that in  Ciuperca(2013), Theorem 1, using also the Lemma \ref{Lemma3 LASSObis} stated in Section 7 and the Remark \ref{Lemma4 Bai}. The details are omitted.
\hspace*{\fill}$\Diamond$ \\

\noindent {\bf Proof of Theorem \ref{Theorem6 LASSOadaptbis}}\\
\textit{(i)} The assertion follows from the  Theorem \ref{Theorem1 LASSObis} and from the Theorem 3(b) in  Xu and Ying (2010).\\
\textit{(ii)} By Xu and Ying (2010), we have: $\lim_{n \rightarrow \infty} \eP \cro{{\cal A}^L_{n;(l^0_{r-1}, l^0_r)}={\cal A}_{(l^0_{r-1},l^0_r)}}=1$. The asymptotic normality of the estimators implies that, for all $k \in {\cal A}_{(l^0_{r-1},l^0_r)} $ we have $\phi^0_{r,k} - \hat \phi^L_{n;(\hat l^L_{r-1}, \hat l^L_r),k} \overset{{\eP}} {\underset{n \rightarrow \infty}{\longrightarrow}} 0$. Thus $k \in {\cal A}^L_{n;(\hat l^L_{r-1}, \hat l^L_r)}$. Hence 
$\lim_{n \rightarrow \infty} \eP \cro{ {\cal A}^L_{n;(\hat l^L_{r-1}, \hat l^L_r)} \supseteq {\cal A}^L_{(l^0_{r-1},l^0_r)}}=1$. 
The proof is finished if we show the claim
$
\eP \cro{\exists k \in \{1, \cdots, p \}, k \not \in {\cal A}^L_{(l^0_{r-1},l^0_r)}, k \in {\cal A}^L_{n;(\hat l^L_{r-1}, \hat l^L_r)} } {\longrightarrow}  0$, as $n \rightarrow \infty$.
Since $k \in {\cal A}^L_{n;(\hat l^L_{r-1}, \hat l^L_r)}$ we have that, with the probability 1,
\begin{equation}
\label{eq10bis}
sgn(\hat \phi^L_{(\hat l^L_{r-1}, \hat l^L_r),k}) \neq 0.
\end{equation}
We suppose, without loss of generality, that $sgn(\hat \phi^L_{(\hat l^L_{r-1}, \hat l^L_r),k}) =1$. Then, using the KKT conditions, we have with the probability 1,
\begin{equation}
\label{eq9bis}
\lambda_{n,(\hat l^L_{r-1}, \hat l^L_r),k} \cdot sgn(\hat \phi^L_{(\hat l^L_{r-1}, \hat l^L_r),k}) =\sum^{\hat l^L_r}_{i=\hat l^L_{r-1}+1} X_{ik} \cdot sgn(Y_i-\XX^t_i \hat \ef^L_{(\hat l^L_{r-1};\hat l^L_r)} ).
\end{equation}
On the other hand, since $k \not \in {\cal A}^L_{(l^0_{r-1},l^0_r)}$, we have $\phi^0_{r,k}=0$, then $sgn(\phi^0_{r,k})=0$. By the proof of the Proposition 4 in the paper of Xu and Ying(2010), for $\phi^0_{,k}=0$, we have that for every  $\ef$ such that $\|\phi_{,k}\| \leq C n^{-1/2}$,
\[
\eP \cro{ sgn \pth{- \sum^n_{i=1} X_{ik} \cdot sgn(Y_i-\XX^t_i \ef)+\lambda_{n,k} \cdot  sgn(\phi_{,k}) } =sgn(\phi_{,k})} {\underset{n \rightarrow \infty}{\longrightarrow}} 1.
\]
Then, taking into account the assertion (i), we apply the previous relation for $\hat \phi^L_{(\hat l^L_{r-1};\hat l^L_r),k}$
\[
\lim_{n \rightarrow \infty} \eP \cro{ sgn \pth{-\sum^n_{i=1} X_{ik} \cdot  sgn(Y_i-\XX^t_i \hat \ef^L_{(\hat l^L_{r-1};\hat l^L_r)})+\lambda_{n,(\hat l^L_{r-1};\hat l^L_r),k} \cdot  sgn(\hat \phi^L_{(\hat l^L_{r-1};\hat l^L_r),k }} =sgn(\phi_{,k})} =1,
\]
where $\hat \phi^L_{(\hat l^L_{r-1};\hat l^L_r),k}$ is the $k$th component of the random vector $\hat \ef^L_{(\hat l^L_{r-1};\hat l^L_r)}$. 
Moreover, by (\ref{eq9bis}), \\ $sgn \pth{-\sum^n_{i=1} X_{ik} \cdot  sgn(Y_i-\XX^t_i \hat \ef^L_{(\hat l^L_{r-1};\hat l^L_r)})+\lambda_{n,(\hat l^L_{r-1};\hat l^L_r),k} \cdot sgn(\hat \phi^L_{(\hat l^L_{r-1};\hat l^L_r),k) })}$ is 0, with the probability 1. Then $\lim_{n \rightarrow \infty} \eP \cro{ 0=sgn(\phi_{,k})}=1$. 
Contradiction with (\ref{eq10bis}). Thus the claim holds.
\hspace*{\fill}$\Diamond$ \\

The demonstration of the  Remark \ref{Lemma4 Bai} is similar to that of Bai(1998). Then we give only the main idea.

\noindent {\bf Proof of Remark \ref{Lemma4 Bai}}\\
Similar as for the proof of the  Lemma \ref{Lemma3 Bai}, we obtain, for $\| \ef_1-\ef_2\|_2 \leq c_n n^{-1/2}$,
\[
\frac{1}{nc_n^2} \left| \pth{{\cal  G}^{(\tau)}_{n}(\ef_1;\ef^0)- {\cal G}^{(\tau)}_{n}(\ef_2;\ef^0)-\eE[{\cal G}^{(\tau)}_{n}(\ef_1;\ef^0)]+\eE[{\cal G}^{(\tau)}_{n}(\ef_2;\ef^0)]} \right| 
\leq C \frac{\| \ef_1-\ef_2\|_2}{c^2_n} \leq \frac{C}{\sqrt{n} c_n},
\]
which converges to 0 for $n \rightarrow \infty$. We also have
$
\eP [ \sup_j | n^{-1} c^{-2}_n  [{\cal G}^{(\tau)}_{n}(\ef_j;\ef^0) -\eE[{\cal G}^{(\tau)}_{n}(\ef_j;\ef^0) ]] | > \epsilon ]
\leq \sum_j \eP [|  [{\cal G}^{(\tau)}_{i}(\ef_j;\ef^0) -\eE[{\cal G}^{(\tau)}_{n}(\ef_j;\ef^0) ]] | > n c^2_n \epsilon]$, where $j=1, \cdots , n^{p/2}$.
By the relations (\ref{eq10}) and (\ref{eq11}), we have
$ |G^{(\tau)}_{i}(\ef;\ef^0) -\eE[G^{(\tau)}_{i}(\ef;\ef^0) ] | \leq C | \XX^t_i (\ef-\ef^0)| < C c_n$. The rest of proof is similar to that of the Lemma 4 of Bai(1998).
\hspace*{\fill}$\Diamond$ \\

\section{Lemmas}
We present in this section the lemmas with proofs, which are useful to prove the main results. 
Following Lemma gives the asymptotic behavior of the objective function $G^{(\tau)}_{i}$ without penalty. In fact, Lemma \ref{Lemma3 Bai} will be necessary to prove the Lemmas \ref{Lemma1 LASSO} and \ref{Lemma1 LASSObis}, where the penalized objective functions are studied.

\begin{lemma}
\label{Lemma3 Bai} Under the assumptions (A1), (A3), 
for all $\alpha >1/2$, we have \\
$
\sup_{1 \leq l < k \leq n} | \inf_{\ef \in \Gamma} \sum^k_{i=l} G^{(\tau)}_{i}(\ef;\ef^0)  | =O_{\eP}(n^\alpha)$.
\end{lemma}
\noindent {\bf Proof of Lemma \ref{Lemma3 Bai}}\\
By direct calculations $R^{(\tau)}_{i}(\ef_1;\ef^0)-R^{(\tau)}_{i}(\ef_2;\ef^0)$ can be written as 
\[
 \acc{\XX^t_i(\ef_1-\ef_2)[(1-\tau)\e1_{\varepsilon_i \leq 0}-\tau \e1_{\varepsilon_i >0}])}  +\acc{ [\varepsilon_i-\XX^t_i(\ef_1-\ef^0)][\tau \e1_{\varepsilon_i > \XX_i^t(\ef_1-\ef^0)}-(1-\tau \e1_{\varepsilon_i \leq \XX^t_i(\ef_1-\ef^0)})  ]} \]
 \begin{equation}
\label{etr}
- \acc{ [\varepsilon_i-\XX^t_i(\ef_2-\ef^0)][\tau \e1_{\varepsilon_i > \XX_i^t(\ef_2-\ef^0)}-(1-\tau \e1_{\varepsilon_i \leq \XX^t_i(\ef_2-\ef^0)})  ] } 
 \equiv S_{1,i} +S_{2,i} - S_{3,i}.
\end{equation}
Obviously $S_{1,i}=\XX^t_i(\ef_1-\ef_2)D_i$. For $S_{2,i} - S_{3,i}$ we have:\\
 \underline{If $\varepsilon_i > \XX_i^t(\ef_1-\ef^0)$}.  When $\varepsilon_i > \XX_i^t(\ef_2-\ef^0)$, we have $S_{2,i} - S_{3,i}=\tau \XX^t_i(\ef_2-\ef_1)$. In the case $\varepsilon_i \leq  \XX_i^t(\ef_2-\ef^0)$, we have $S_{2,i} - S_{3,i} = \tau \XX^t_i(\ef_2-\ef_1) +[\varepsilon_i-\XX^t_i(\ef_2-\ef^0)] \leq \tau \XX^t_i(\ef_2-\ef_1)$.
Then, in the both cases, $S_{2,i} - S_{3,i} \leq \tau \XX^t_i(\ef_2-\ef_1)$.\\
\underline{If $\varepsilon_i \leq \XX_i^t(\ef_1-\ef^0)$}.  When $\varepsilon_i \leq  \XX_i^t(\ef_2-\ef^0)$, we have $S_{2,i} - S_{3,i}=(1-\tau) \XX^t_i(\ef_1-\ef_2)$. In the case $\varepsilon_i >  \XX_i^t(\ef_2-\ef^0)$, we have $S_{2,i} - S_{3,i} = (1-\tau)[ \XX^t_i(\ef_1-\ef_2) +\varepsilon_i-\XX^t_i(\ef_1-\ef^0)] \leq (1-\tau) \XX^t_i(\ef_1-\ef_2)$.
Then, in the both cases, $S_{2,i} - S_{3,i} \leq (1-\tau) \XX^t_i(\ef_1-\ef_2)$.\\
In conclusion, with the probability 1,
\begin{equation}
\label{eq10}
S_{1,i} +S_{2,i} - S_{3,i} \leq \XX^t_i(\ef_1-\ef_2)D_i +\max (\tau \XX^t_i(\ef_2-\ef_1),(1-\tau)\XX^t_i(\ef_1-\ef_2))
\end{equation}
Similarly
\begin{equation}
\label{eq11}
S_{1,i} +S_{2,i} - S_{3,i} \geq \XX^t_i(\ef_1-\ef_2)D_i +\min (\tau \XX^t_i(\ef_2-\ef_1),(1-\tau)\XX^t_i(\ef_1-\ef_2))
\end{equation}
Hence, the relations (\ref{etr}), (\ref{eq10}) and (\ref{eq11}), for $\|\ef_1-\ef_2 \|_2 \leq C n^{-1/2}$ together the assumption (A3), imply that
\begin{equation}
\label{eq12}
\begin{array}{c}
\left|\sum^n_{i=1}  \cro{R^{(\tau)}_{i}(\ef_1;\ef^0)-R^{(\tau)}_{i}(\ef_2;\ef^0)-\eE[R^{(\tau)}_{i}(\ef_1;\ef^0)]+\eE[R^{(\tau)}_{i}(\ef_2;\ef^0)] }  \right|  \\
\leq C \sum^n_{i=1} \|\XX_i\|_2 \cdot \| \ef_2-\ef_1\|_2 \leq O_{\eP}(n^{1/2}).
\end{array}
\end{equation}
By an argument similar to the one used in the Lemma 3 of Bai(1998),  together the Proposition \ref{inegE0},    we obtain
\[
\sup_{1 \leq l < k \leq n} \left| \inf_\ef \sum^k_{i=l} G^{(\tau)}_{i}(\ef;\ef^0) \right| \leq 2 \sup_{1 \leq k \leq n} \sup_\ef \left| \sum^k_{i=1} [ G^{(\tau)}_{i}(\ef;\ef^0) - \eE[G^{(\tau)}_{i}(\ef;\ef^0)]] \right |.
\]
On the other hand
$
\sum^k_{i=1} [ G^{(\tau)}_{i}(\ef;\ef^0) - \eE[G^{(\tau)}_{i}(\ef;\ef^0)]]=\sum^k_{i=1} [ R^{(\tau)}_{i}(\ef;\ef^0) - \eE[R^{(\tau)}_{i}(\ef;\ef^0)]]+\sum^k_{i=1} D_i \XX^t_i(\ef-\ef^0).
$ 
Let us consider the random process
$
\xi_k=\sup_\ef \left|\sum^k_{i=1} [ G^{(\tau)}_{i}(\ef;\ef^0) - \eE[G^{(\tau;\ef^0)}_{i}(\ef;\ef^0)]]  \right|
$. 
Then, since by Proposition \ref{inegE0} $\eE[G^{(\tau)}_{i}(\ef;\ef^0)] \geq 0$,  follows that $\acc{\xi_k, {\cal F}_k}_{k=1, \cdots, n}$ is a sub-martingale, where ${\cal F}_k=\sigma-field\{\varepsilon_1, \cdots, \varepsilon_k\} $, which implies, using Doob's inequality 
$
\eP[ \sup_{1 \leq k \leq n} \xi_k > n^\alpha] \leq n^{-\alpha m} C_m \eE[\xi^m_n], \qquad C_m >0$,
with $m>1$.
We divide the parameter set $\Gamma$ into $m^{p/2}$ cells, such that the cell diameter is  $\leq n^{-1/2}$. Thus
$ 
| \sum^n_{i=1} [G^{(\tau)}_{i}(\ef_1;\ef^0) - \eE[G^{(\tau)}_{i}(\ef_1;\ef^0)] - G^{(\tau)}_{i}(\ef_2;\ef^0) + \eE[G^{(\tau)}_{i}(\ef_2;\ef^0)]] |
$
$
\leq | \sum^n_{i=1} \cro{R^{(\tau)}_{i}(\ef_1;\ef^0) - \eE[R^{(\tau)}_{i}(\ef_1;\ef^0)] - R^{(\tau)}_{i}(\ef_2;\ef^0) + \eE[R^{(\tau)}_{i}(\ef_2,\ef^0)]} | +| \sum^n_{i=1} D_i \XX^t_i (\ef_1 - \ef_2) |
$
and using the  relation (\ref{eq12}), we obtain, with the probability 1,  that the last relation is smaller than
$ \sum^n_{i=1} | \XX^t_i(\ef_2-\ef_1)| \leq C n n^{-1/2}=C n^{1/2}$.
By an argument similar to the one used in  Bai(1998), we have $\eE \left| \sum^n_{i=1} [G^{(\tau)}_{i}(\ef_r; \ef^0) -\eE[G^{(\tau)}_{i}(\ef_r;\ef^0)]] \right| \leq C n^{m/2}$. The rest of proof is similar to that of the Lemma 3 of Bai(1998).
\hspace*{\fill}$\Diamond$ \\

\subsection{For SCAD estimator}
Following result will be useful in the study of the convergence rate of the change-point SCAD estimator in a model with breaking. 
\begin{lemma}
\label{Lemma1 LASSO}
Under the assumptions (A1), (A3), for a positive sequence $(\lambda_n)_n$ such that  $\lambda_n \rightarrow 0$, we have 
\[
\sup_{0 \leq j_1 <j_2 \leq n} \left| \inf_\ef \sum^{j_2}_{i=j_1+1} G^{(\tau,\lambda)}_{i} (\ef;\ef^0)\right|=O_{\eP}(n^\alpha,n \lambda_n).
\]
\end{lemma}
\noindent {\bf Proof of Lemma \ref{Lemma1 LASSO}}\\
Using the triangle inequality, we deduce that
\[
\sup_{0 \leq j_1 <j_2 \leq n} \left| \inf_\ef \sum^{j_2}_{i=j_1+1} G^{(\tau,\lambda)}_{i} (\ef;\ef^0)\right| \leq \sup_{0 \leq j_1 <j_2 \leq n} \left| \inf_\ef \sum^{j_2}_{i=j_1+1} G^{(\tau)}_{i} (\ef;\ef^0)\right| + n\sup_\ef \left| (\pp_\lambda(|\ef|)- \pp_\lambda(| \ef^0|))^t \uu_p  \right|. 
\]
Considering  Lemma \ref{Lemma3 Bai} and the definition of $p'_\lambda$, we have that the last quantity is smaller than  $O_{\eP}(n^\alpha) +n \lambda_n$.
\hspace*{\fill}$\Diamond$ \\

In the following Lemma, the behavior of ${\cal G}^{(\tau)}_{n}$ is studied in the outside of the ball center $\eb^0$ and radius $c_n$. 
\begin{lemma}
\label{Lemma5 Bai}
Under the assumptions (A1), (A2), with $(c_n)$ a positive sequence  such that  $c_n \rightarrow 0$ and $n c^2_n/ \log n\rightarrow \infty$, there exists $\epsilon > 0$ such that we have with probability 1
\[
\liminf_{n \rightarrow \infty} \pth{\inf_{\|\ef-\ef^0 \|_2 \geq  c_n } \frac{1}{n c^2_n}{\cal G}^{(\tau)}_{n}(\ef;\ef^0) } \geq \epsilon > 0.
\]
\end{lemma}
\noindent {\bf Proof of Lemma \ref{Lemma5 Bai}}\\
Let $\UU$ in an open subset of $\R^p$.
By the proof of the  Lemma 3  of Wu and Liu(2009), taking into account the assumptions (A1) and (A2),  we have
\[
\eE\cro{{\cal G}^{(\tau)}_{n} \pth{ \ef^0+\frac{\UU}{\sqrt{n}}; \ef^0} }=\frac{f(0)}{2 n}\UU^t (\sum^n_{i=1} \XX_i \XX_i^t)\UU+o(1).
\]
If $c_n \rightarrow 0$ and $n c^2_n \rightarrow \infty$, we have similarly $\eE\cro{{\cal G}^{(\tau)}_{n} \pth{ \ef^0+\UU c_n  ;\ef^0}}=\frac{f(0)}{2 } c^2_n\UU^t (\sum^n_{i=1} \XX_i \XX_i^t)\UU+o_{\eP}(1).$
The function $G^{(\tau)}_{i}(\ef;\ef^0)$ is convex, hence ${\cal G}^{(\tau)}_{n}(\ef;\ef^0)$ is convex  in $\ef$. Thus, its minimum over  $\| \ef-\ef^0\|_2 \geq c_n$ is realized for  $\| \ef-\ef^0\|_2 = c_n$. Then, for $\| \UU\|_2=1$, using the assumption (A2) we obtain that $\eE\cro{{\cal G}^{(\tau)}_{n} \pth{ \ef^0+\UU c_n ;\ef^0} }=\frac{f(0)}{2} n c^2_n (C+o(1))$. The rest of proof follows using the the Lemma 5 of Bai(1998), taking into account the relation (\ref{as}).
\hspace*{\fill}$\Diamond$ \\
\begin{lemma}
\label{Lemma2 LASSO}
Under the assumptions (A1), (A2), for two positive sequences  $(c_n)$ and $(\lambda_n)$ such that  $\lambda_n \rightarrow 0$, $c_n \rightarrow 0$, $n c^2_n/ \log n\rightarrow \infty$ and $\lambda_n c^{-2}_n  \rightarrow 0$, we have
\[
\liminf_{n \rightarrow \infty} \pth{ \inf_{\|\ef-\ef^0\|_2 \geq c_n} \frac{1}{n c^2_n} \sum^n_{i=1} G^{(\tau,\lambda)}_{i}(\ef;\ef^0) } > \epsilon.
\]
\end{lemma}
\noindent {\bf Proof of Lemma \ref{Lemma2 LASSO}}\\
Applying the mean value theorem, we write
$G^{(\tau,\lambda)}_{i}(\ef;\ef^0)=G^{(\tau)}_{i}(\ef;\ef^0)+ [|\ef|-|\ef^0| ]^t \pp'_\lambda(\tilde \ef)$, with $\tilde \ef =\ef^0-\textbf{b}(\ef-\ef^0)$, $\textbf{b} \in [0,1]^p$. Then, using the relation (\ref{eq14}), we have
\[
\inf_{\|\ef-\ef^0\|_2 \geq c_n}\frac{1}{n c^2_n} \sum^n_{i=1} { G}^{(\tau,\lambda)}_{i}(\ef;\ef^0) 
\geq \inf_{\|\ef-\ef^0\|_2 \geq c_n}\frac{1}{n c^2_n} {\cal G}^{(\tau)}_{n}(\ef;\ef^0)- \frac{n \lambda_n}{n c^2_n}.
\]
Since $\lambda_n c^{-2}_n  \rightarrow 0$, for every $\epsilon >0$ there exists a  $n_\epsilon \in \N$ such that $\lambda_n c^{-2}_n  < \epsilon/2$. An application of Lemma \ref{Lemma5 Bai} leads to
$
\inf_{\|\ef-\ef^0\|_2 \geq c_n}\frac{1}{n c^2_n} {\cal G}^{(\tau)}_{n}(\ef;\ef^0) > \frac{3 \epsilon}{2}
$
and  Lemma is proved.
\hspace*{\fill}$\Diamond$ \\

By similar calculus as in  Bai(1998), Lemma 10, we have following result for the estimator (\ref{eq3}) of $\ef$.

\begin{lemma}
\label{Lemma3 LASSO}
For $(\lambda_n)$, $(c_n)$ as in the Lemma \ref{Lemma2 LASSO}), under assumptions (A1)-(A3), for all $n_1, n_2 \in \N$ such that $n_1 \geq n^u$, with $3/4 \leq u  \leq 1$, $ n_2 \leq n^v$,  $v < 1/4$, let us consider the model
\[
\begin{array}{lll}
Y_i=\XX_i^t \ef^0_1+\varepsilon_i, & &i=1, \cdots, n_1 \\
Y_i=\XX_i^t \ef^0_2+\varepsilon_i, & &i=n_1+1, \cdots, n_1+n_2 \\
\end{array}
\] 
with the assumption $\ef^0_1 \neq \ef^0_2$.  
Consider $A^{(\tau,\lambda)}_{n_1+n_2}(\ef)=\sum^{n_1}_{i=1} G^{(\tau,\lambda)}_{i} (\ef;\ef^0_1)+ \sum^{n_1+n_2}_{i=n_1+1} G^{(\tau,\lambda)}_{i}(\ef;\ef^0_2)$ and
$\hat \ef^{(\tau,\lambda)}_{n_1+n_2} \equiv \argmin _{\ef} A^{(\tau,\lambda)}_{n_1+n_2}(\ef)$.\\
(i) $\|\hat \ef^{(\tau,\lambda)}_{n_1+n_2} - \ef^0_1 \|_2  \leq n_1^{-1/2} n_1^{\frac{v+\delta}{2u}} \leq n^{-(u-v-\delta)/2}$. \\
(ii) $\sum^{n_1}_{i=1}G^{(\tau,\lambda)}_{i} \pth{\hat \ef^{(\tau,\lambda)}_{n_1+n_2}; \ef^0_1 }=O_{\eP}(1)$.
\end{lemma}
\noindent {\bf Proof of Lemma \ref{Lemma3 LASSO}}\\
\textit{(i)} $A^{(\tau,\lambda)}_{n_1+n_2} (\hat \ef^{\tau,\lambda}_{n_1+n_2} ) \leq \sum^{n_1+n_2}_{i=n_1+1} G^{(\tau,\lambda)}_{i}(\ef^0_1;\ef^0_2)=\sum^{n_1+n_2}_{i=n_1+1} G^{(\tau)}_{i}(\ef^0_1;\ef^0_2)+n_2 [ \pp_\lambda(|\ef^0_1|)-\pp_\lambda(|\ef^0_2|)]\uu_p=o_{\eP}(1)+O(n_2)$. By Lemma \ref{Lemma2 LASSO}, for $G^{(\tau,\lambda)}_{i}$, $i=n_1+1, \cdots, n_2$,  we arrive to a contradiction.\\
\textit{(ii)} Let $Z^{(\tau)}_n(\ef)\equiv \sum^{n_1}_{i=1} G^{(\tau)}_{i}(\ef;\ef^0_1)$, $t^{(\tau)}_n(\ef) \equiv\sum^{n_1+n_2}_{i=n_1+1} \cro{\rho_\tau(\varepsilon_i-\XX^t_i (\ef-\ef^0_2)) -\rho_\tau(\varepsilon_i-\XX^t_i(\ef^0_1-\ef^0_2)) }$, $t^{(\tau,\lambda)}_n(\ef) \equiv t^{(\tau)}_n(\ef) +n_2 [\pp_\lambda(|\ef|)- \pp_\lambda(|\ef^0_1|)] \uu_p $, $Z^{(\tau,\lambda)}_n(\ef)\equiv Z^{(\tau)}_n(\ef)+n_1 [\pp_\lambda(|\ef|)- \pp_\lambda(|\ef^0_1|)] \uu_p$. Then
\[
A^{(\tau,\lambda)}_{n_1+n_2}(\ef) =Z^{(\tau,\lambda)}_n(\ef) +t^{(\tau,\lambda)}_n(\ef) +n_2 [\pp_\lambda(|\ef^0_1|)- \pp_\lambda(|\ef^0_2|)] \uu_p +\sum^{n_1+n_2}_{i=n_1+1} \cro{\rho_\tau(\varepsilon_i-\XX^t_i (\ef^0_1-\ef^0_2)) -\rho_\tau(\varepsilon_i) }.
\]
We have
$| t^{(\tau,\lambda)}_n(\hat \ef^{(\tau,\lambda)}_{n_1+n_2})| \leq |  t^{(\tau)}_n(\hat \ef^{(\tau,\lambda)}_{n_1+n_2})|+n_2  \| \pp_\lambda(| \hat \ef^{(\tau,\lambda)}_{n_1+n_2}|) - \pp_\lambda(|\ef^0_1 |) \|_1$ and similarly that for the relation (\ref{eq10}), $\leq C \sum^{n_1+n_2}_{i=n_1+1} |\XX_i|_1 \| \hat \ef^{(\tau,\lambda)}_{n_1+n_2} -\ef^0_1 \|_2+n_2 \| \pp_\lambda(| \hat \ef^{(\tau,\lambda)}_{n_1+n_2}|) - \pp_\lambda(|\ef^0_1 |) \|_1$. The rest of proof is similar to that of the  Lemma 3(ii) of Ciuperca(2013), taking into account the assumption (A3).
\hspace*{\fill}$\Diamond$ \\

We have the equivalent of  Lemma 4 of the same paper.\\

\subsection{For LASSO-type estimator}
\begin{lemma}
\label{Lemma1 LASSObis}
Under the assumptions (A1), (A3), we have, for $\alpha > 1/2$,
\[
\sup_{0 \leq j_1 < j_2 \leq n} \left| \inf_{\ef} \sum^{j_2}_{i=j_1+1} \eta^L_{i;(j_1,j_2)}(\ef, \ef^0) \right|= O_{\eP}(\max(n^\alpha, \sup_{0 \leq j_1 < j_2 \leq n} \|\el_{n;(j_1,j_2)}\|_2)).
\] 
\end{lemma}
\noindent {\bf Proof of Lemma \ref{Lemma1 LASSObis}}\\
By the  Lemma 3 of Bai(1998),  Lemma holds for $G^{(1/2)}_{i}$ instead of $ \eta^L_{i;(j_1,j_2)}$. For $\eta^L_i$, we have $| \eta^L_{i;(j_1,j_2)}(\ef; \ef^0) - G^{(1/2)}_{i}(\ef;\ef^0)| =(j_1-j_1)^{-1/2} | \el_{n;(j_1,j_2)}^t (|\ef|-|\ef^0|)|$. 
Then, by triangular inequality together the Lemma \ref{Lemma3 Bai} and the compactness of the set $\Gamma$
\[
\sup_{0 \leq j_1 < j_2 \leq n} \left| \inf_{\ef} \sum^{j_2}_{i=j_1+1} \eta^L_{i;(j_1,j_2)}(\ef; \ef^0) \right| \leq  \sup_{0 \leq j_1 < j_2 \leq n} \left| \inf_{\ef} \sum^{j_2}_{i=j_1+1} G_i^{(1/2)}(\ef; \ef^0) \right| +C \sup_{0 \leq j_1 < j_2 \leq n} \| \el_{n;(j_1,j_2)}\|_1
\]
$ =O_{\eP}(n^\alpha)+O_{\eP}(\el_n)$.
\hspace*{\fill}$\Diamond$ \\
\begin{lemma}
\label{Lemma2 LASSObis}
Under the assumptions (A1), (A2), if $n^{-1/2}\| \el_n\|_2 \overset{{\eP}} {\underset{n \rightarrow \infty}{\longrightarrow}} \lambda_0 $, with $\lambda_0 \geq 0$, then
\[
\liminf_{n \rightarrow \infty} \pth{\inf_{\| \ef-\ef^0 \|_2 \geq n^{-1/2}} n^{-1} \sum^n_{i=1} \eta^L_{i;(0,n)} (\ef; \ef^0)} > \epsilon.
\]
\end{lemma}
\noindent {\bf Proof of Lemma \ref{Lemma2 LASSObis}}\\
Using (A2), by the Lemma 6 of Ciuperca(2011b), we have for $G^{(1/2)}_{i}$, with the probability 1
\begin{equation}
\label{eq27}
\liminf_{n \rightarrow \infty} \pth{\inf_{\| \ef-\ef^0 \|_2 \geq n^{-1/2}} n^{-1} \sum^n_{i=1} G^{(1/2)}_{i} (\ef; \ef^0)} > \frac{3 \epsilon}{2}.
\end{equation}
We also have the inequality
\begin{equation}
\label{eq28}
\inf_{\| \ef-\ef^0 \|_2 \geq n^{-1/2}}  \sum^n_{i=1} \eta^L_{i;(0,n)} (\ef; \ef^0) \geq \inf_{\| \ef-\ef^0 \|_2 \geq n^{-1/2}}  \sum^n_{i=1} G^{(1/2)}_{i} (\ef; \ef^0) - \sup_{\| \ef-\ef^0 \|_2 \geq n^{-1/2}} \el^t_n  ( |\ef| -|\ef^0|).
\end{equation}
Since $\| \el_n\|_2=O_{\eP}(n^{-1/2})$ and $\ef $ belongs to a  compact set, then the last term of the right-hand-side of (\ref{eq28}) is $o_{\eP}(n^{-1})$. Hence
\begin{equation}
\label{eq29}
\sup_{\| \ef-\ef^0 \|_2 \geq n^{-1/2}} \pth{ n^{-1} \el^t_n( |\ef| -|\ef^0|) } <\frac{\epsilon}{2}, \qquad n \geq n_\epsilon .
\end{equation}
The conclusion follows,  combining the relations (\ref{eq27}), (\ref{eq28}) and (\ref{eq29}).  
\hspace*{\fill}$\Diamond$ \\
\begin{lemma}
\label{Lemma3 LASSObis}
For all $n_1, n_2 \in \N$ such that $n_1 \geq n^u$, with $3/4 \leq u  \leq 1$, $ n_2 \leq n^v$,  $v < 1/4$, let us consider the model
\[
\begin{array}{lll}
Y_i=\XX_i^t \ef^0_1+\varepsilon_i, & &i=1, \cdots, n_1 \\
Y_i=\XX_i^t \ef^0_2+\varepsilon_i, & &i=n_1+1, \cdots, n_1+n_2 \\
\end{array}
\] 
with the assumption $\ef^0_1 \neq \ef^0_2$.  
Under the assumptions (A1)-(A3) and $(\lambda_n)$ as in the Lemma \ref{Lemma2 LASSObis}, let us consider
$
A^L_{n_1+n_2}(\ef)=\sum^{n_1}_{i=1} \eta^L_{i;(0,n_1)} (\ef; \ef^0_1)+\sum^{n_1+n_2}_{i=n_1+1} \eta^L_{i;(n_1,n_1+n_2)} (\ef; \ef^0_2)$
 and $\hat \ef^L_{n_1+n_2} \equiv \argmin_{\ef} A^L_{n_1+n_2}(\ef)$. Then\\
 (i) $\| \hat \ef^L_{n_1+n_2}- \ef^0_1 \|_2 \leq  n_1^{-1/2} n_1^{\frac{v+\delta}{2u}} \leq n^{-(u-v-\delta)/2}$.  \\
 (ii) $\sum^{n_1}_{i=1} \eta^L_{i;(0,n_1)}( \hat \ef^L_{n_1+n_2}, \ef^0_1)=O_{\eP}(1)  $.
\end{lemma}
\noindent {\bf Proof of Lemma \ref{Lemma3 LASSObis}}\\
We denote  $\hat \ef_{n_1+n_2} \equiv \argmin_{\ef} \sum^{n_1}_{i=1} G^{(1/2)}_{i} (\ef, \ef^0_1)+\sum^{n_1+n_2}_{i=n_1+1} G^{(1/2)}_{i} (\ef, \ef^0_2)$. Using the assumptions (A1) and (A3), by Lemma 10 of Bai(1998) we have that  (i) and (ii) are true for  $\hat \ef_{n_1+n_2}$ and $G^{(1/2)}_{i}$.\\
(i) We suppose the contrary  $\| \hat \ef^L_{n_1+n_2}- \ef^0_1 \|_2 \geq n_1$. On the other hand, we have by definition:
\begin{equation}
\label{eq31}
A^L_{n_1+n_2}(\hat \ef^L_{n_1+n_2}) \leq \sum^{n_1+n_2}_{i=n_1+1} \eta^L_{i;(n_1,n_1+n_2)} (\ef^0_1, \ef^0_2) =\sum^{n_1+n_2}_{i=n_1+1} G^{(1/2)}_{i} (\ef^0_1, \ef^0_2) +\el^t_n( |\ef^0_1| -|\ef^0|) .
\end{equation}
By Lemma 10(ii) of Bai(1998), we have: $\sum^{n_1+n_2}_{i=n_1+1} G^{(1/2)}_{i} (\ef^0_1, \ef^0_2) =o_{\eP}(1)$, then taking into account the relation  (\ref{eq31}), we obtain
\begin{equation}
\label{eq32}
A^L_{n_1+n_2}(\hat \ef^L_{n_1+n_2}) \leq o_{\eP}(1)+O_{\eP}(n^{1/2}).
\end{equation}
On the other hand, using the  Lemma \ref{Lemma2 LASSObis}, we deduce
\begin{equation}
\label{eq33}
\sum^{n_1}_{i=1} \eta^L_{i;(0,n_1)}( \hat \ef^L_{n_1+n_2}, \ef^0_1) \geq O_{\eP}(n_1).
\end{equation}
There is a contradiction between the relations (\ref{eq32}) and (\ref{eq33}).\\
(ii) Introduce $\nu_{(\ef_1,\ef_2)}(\XX_i) \equiv \XX^t_i(\ef_1-\ef_2)$. For $\tau=1/2$, let us recall the  notations given in  Lemma \ref{Lemma3 LASSO}: $Z^{(1/2)}_n(\ef) \equiv \sum^{n_1}_{i=1} G^{(1/2)}_{i}(\ef, \ef^0_1)$,
 $t^{(1/2)}_n(\ef) \equiv \sum^{n_1+n_2}_{i=n_1+1} [ |\varepsilon_i- \nu_{(\ef,\ef^0_2)}(\XX_i) |-| \varepsilon_i- \nu_{(\ef^0_1,\ef^0_2)}(\XX_i)|]$. By the Lemma 7 of Ciuperca(2011b), we have that: $Z^{(1/2)}_n(\hat \ef^L_{n_1+n_2},\ef^0_1)=O_{\eP}(1)$. Introduce now  $t^L_n(\ef) \equiv t^{(1/2)}_n(\ef)+\el^t_{(n_1,n_1+n_2)}[|\ef|-|\ef^0_1|]$, $Z^L_n(\ef) \equiv Z^{(1/2)}_n(\ef)+\el^t_{(0,n_1)}[|\ef|-|\ef^0_1|]$. Thus
 $
 A^L_{n_1+n_2}(\ef)=Z^{(1/2)}_n(\ef)+\el^t_{(0,n_1)}[|\ef|-|\ef^0_1|] +t^{(1/2)}_n(\ef)+\sum^{n_1+n_2}_{i=n_1+1}  \cro{|\varepsilon_i-\nu_{(\ef^0_1,\ef^0_2)}(\XX_i) | -|\varepsilon_i |} 
 +\el^t_{(n_1,n_1+n_2)}[|\ef|-|\ef^0_2|]
 $ $
 =Z_n^L(\ef) +t^L_n(\ef)+\el^t_{(n_1,n_1+n_2)}[|\ef^0_1|-|\ef^0_2|] - \sum^{n_1+n_2}_{i=n_1+1}  \cro{|\varepsilon_i-\nu_{(\ef^0_1,\ef^0_2)}(\XX_i) | -|\varepsilon_i |}$. 
 Then $
 \hat \ef^L_{n_1+n_2} \equiv \argmin_{\ef} A^L_{n_1+n_2}(\ef)=\argmin_{\ef} [Z_n^L(\ef) +t^L_n(\ef) ]$. 
 But $| t^L_n( \hat \ef^L_{n_1+n_2}) | \leq | t^{(1/2)}_n( \hat \ef^L_{n_1+n_2}) | +| \el^t_{(n_1,n_1+n_2)}[|\hat \ef^L_{n_1+n_2}|-|\ef^0_1|]| $ and using the elementary  inequality $||a|-|b| | \leq |a-b|$, we have $ | t^L_n( \hat \ef^L_{n_1+n_2}) | \leq \sum^{n_1+n_2}_{i=n_1+1}\| \hat \ef^L_{n_1+n_2} - \ef^0_1 \|_1 \cdot \|\XX^t_i\|_1 +| \el^t_{(n_1,n_1+n_2)}[|\hat \ef^L_{n_1+n_2}|-|\ef^0_1|]| \leq o_{\eP}(1)$, we have used  (i) and the assumptions (A3) and $\| \el_{(n_1,n_1+n_2)} \|_2=o_{\eP}(n_2^{1/2})$.\\
 We have also $Z_n^L(\ef^0_1) =t^L_n(\ef^0_1)$. Thus $0 \geq \inf_{\ef}( Z^L_n(\ef)+t^L_n(\ef)  )= Z^L_n( \hat \ef^L_{n_1+n_2})+ t^L_n( \hat \ef^L_{n_1+n_2}) = Z^L_n( \hat \ef^L_{n_1+n_2}) -|o_{\eP}(1) | \geq \inf_{\ef} Z^L_n(\ef)- |o_{\eP}(1) | $. But $\inf_{\ef} Z_n(\ef)=O_{\eP}(1)$. The rest of proof is similar to that of the  Lemma 3(ii) of Ciuperca(2011a).
\hspace*{\fill}$\Diamond$ \\

\vskip 5cm
\begin{center}

\textbf{Appendix}
\end{center}
\begin{table}[!p]
  \caption{\footnotesize Median of change-point estimations, percentage of true 0 and of false 0 by LS, QUANT, QLASSO, SCAD,   LASSO-type and  adaptive LASSO methods for $n=200$, $K=2$, $l^0_1=30$, $l^0_2=100$, $\varepsilon_i \sim {\cal E}xp(-1.5,1)$.}
\begin{center}
\begin{tabular}{|c|c|c|c|c|c|c|} \hline 
Method & LS  & QUANT & QLASSO &SCAD  & LASSO-type & aLASSO \\ \hline
   median of $(\hat l_1,\hat l_2)$  & (31,100) & (31,100)   &(31,100) & (30,100) & (30,100) & (30,100)  \\ \hline
   $\%$ of trues 0 & 0 & 0& 46  &  75 & 97& 94\\ \hline
   $\%$ of   false 0 & 0& 0 & 1 & 3 & 3 &  8   \\ \hline
\end{tabular} 
\end{center}
\label{Tabl1} 
\end{table}

\begin{table}[!p]
 \caption{\footnotesize Median of change-point estimations, percentage of true 0 and of false 0 by LS, QUANT, QLASSO, SCAD,   LASSO-type and  adaptive LASSO methods for $n=200$, $K=2$, $l^0_1=30$, $l^0_2=100$, $\varepsilon_i \sim {\cal N}(0,1)$.}
\begin{center}
\begin{tabular}{|c|c|c|c|c|c|c|} \hline 
Method & LS  & QUANT & QLASSO &SCAD  & LASSO-type & aLASSO \\ \hline
   median of $(\hat l_1,\hat l_2)$  & (31,100) & (30,100)   &(30,100) & (30,100) & (30,100) & (30,100)  \\ \hline
   $\%$ of trues 0 & 0 & 0& 37 &  65 & 98& 94\\ \hline
   $\%$ of   false 0 & 0& 0 & 0.5 & 5 & 2 &  8   \\ \hline
\end{tabular} 
\end{center}
\label{Tabl2} 
\end{table}
\begin{table}[!p]
 \caption{\footnotesize Median of change-point estimations, percentage of true 0 and of false 0 by LS, QUANT, QLASSO, SCAD,   LASSO-type and  adaptive LASSO methods for $n=200$, $K=2$, $l^0_1=30$, $l^0_2=100$, $\varepsilon_i \sim Cauchy$.}
\begin{center}
\begin{tabular}{|c|c|c|c|c|c|c|} \hline 
Method & LS  & QUANT & QLASSO &SCAD  & LASSO-type & aLASSO \\ \hline
   median of $(\hat l_1,\hat l_2)$ & (31,100) & (30.5,100)   &(30,100) & (30,100) & (30,100) & (30,100)  \\ \hline
   $\%$ of trues 0 & 0 & 0& 36 &  62 & 95& 48\\ \hline
   $\%$ of   false 0 & 0& 0 & 1 & 3 & 3 &  12   \\ \hline
\end{tabular} 
\end{center}
\label{Tabl3} 
\end{table}
\begin{table}
 \caption{\footnotesize Median of change-point estimations, percentage of true 0 and of false 0 by LS, QUANT, QLASSO,  LASSO-type and  adaptive LASSO methods for $n=60$, $K=2$, $l^0_1=17$, $l^0_2=40$, $\varepsilon_i \sim {\cal E}xp(-1.5,1)$.}
\begin{center}
\begin{tabular}{|c|c|c|c|c|c|} \hline 
Method & LS  & QUANT & QLASSO  & LASSO-type & aLASSO \\ \hline
   median of$(\hat l_1,\hat l_2)$  & (18,40) & (18,40)  &(18,40)  & (18,40) & (17,40)  \\ \hline
   $\%$ of trues 0 & 0 & 0 & 60 &  91 & 75\\ \hline
   $\%$ of   false 0 & 0& 0 & 27 & 27 & 17   \\ \hline
\end{tabular} 
\end{center}
\label{Tabl4} 
\end{table}
\begin{table}
 \caption{\footnotesize Median of change-point estimations, percentage of true 0 and of false 0 by LS, QUANT, QLASSO, LASSO-type and  adaptive LASSO methods for $n=60$, $K=2$, $l^0_1=17$, $l^0_2=40$, $\varepsilon_i \sim {\cal N}(0,1)$.}
\begin{center}
\begin{tabular}{|c|c|c|c|c|c|} \hline 
Method & LS  & QUANT & QLASSO  & LASSO-type & aLASSO \\ \hline
   median of$(\hat l_1,\hat l_2)$  & (18,40) & (18,40)  &(18,40)  & (17,40) & (17,40)  \\ \hline
   $\%$ of trues 0 & 0 & 0 & 51 &  92.5 & 82\\ \hline
   $\%$ of   false 0 & 0& 0 & 7 & 13 & 15   \\ \hline
\end{tabular} 
\end{center}
\label{Tabl5} 
\end{table}
\begin{table}
 \caption{\footnotesize Median of change-point estimations, percentage of true 0 and of false 0 by LS, QUANT, QLASSO,  LASSO-type and  adaptive LASSO methods for $n=60$, $K=2$, $l^0_1=17$, $l^0_2=40$, $\varepsilon_i \sim Cauchy$.}
\begin{center}
\begin{tabular}{|c|c|c|c|c|c|} \hline 
Method & LS  & QUANT & QLASSO  & LASSO-type & aLASSO \\ \hline
   median of$(\hat l_1,\hat l_2)$  & (17,40) & (18,40)  &(18,40)  & (17,40) & (17,40)  \\ \hline
   $\%$ of trues 0 & 0 & 0 & 48.5 &  82 & 43 \\ \hline
   $\%$ of   false 0 & 0& 0 & 18 & 26 & 16   \\ \hline
\end{tabular} 
\end{center}
\label{Tabl6} 
\end{table}

\begin{table}
 \caption{\footnotesize The average of estimation error $\| \hat \ef- \ef^0\|_1$ in each segment under different distributions for LASSO-type and adaptive LASSO methods, $n=200$, $K=2$, $l^0_1=30$, $l^0_2=100$.  }
\begin{center}
\begin{tabular}{|c||ccc|ccc|}\hline
  & LASSO- &type  &  &adaptive   & LASSO &    \\ \hline  
 & $(1,l_1)$  & $(l_1;l_2)$ & $(l_2,n)$  & $(1,l_1)$  & $(l_1;l_2)$ & $(l_2,n)$ \\ \hline 
$\varepsilon_i \sim {\cal N}(0,1)$ & 0.33 & 0.11 & 0.12 & 0.74 & 0.48 & 0.47  \\
$\varepsilon_i \sim {\cal E}xp(-1.5,1)$ & 0.38& 0.11 & 0.12 & 0.82 & 0.52 & 0.48 \\
$\varepsilon_i \sim Cauchy$ & 0.51 & 0.17 & 0.13 & 4.6 & 4.84 & 4.83 \\ \hline
\end{tabular} 
\end{center}
\label{Tabl7} 
\end{table}


\begin{thebibliography}{3}
\bibitem{Babu:89}
 G.J. Babu, (1989), 
\newblock  Strong representations for LAD estimators in linear models.
\newblock {\em Probability Theory and Related Fields}, {\bf 83},  pp.  547-558.
\bibitem{Bai:98}
Bai, J.,(1998). 
\newblock Estimation of multiple-regime regressions with least absolute deviation.
\newblock{\it Journal of Statistical Planning Inference}, \textbf{74},  103-134.
\bibitem{Ciuperca:13}
Ciuperca, G.,(2013). 
\newblock Model selection by LASSO methods in a change-point model, 
\newblock{\it Statistical Papers}, DOI 10.1007/s00362-012-0482-x.
\bibitem{Ciuperca:11a}
 Ciuperca, G.,(2011a). 
\newblock Penalized least absolute deviations estimation for  nonlinear model with change-points.
\newblock {\it Statistical Papers},  \textbf{52}(2), 371-390.
\bibitem{Ciuperca:11b}
Ciuperca, G.,(2011b). 
\newblock Estimating nonlinear regression with and without change-points by the LAD-method.
\newblock{\it Annals of the Institute of Statistical Mathematics},   \textbf{63}(4), 717-743.
\bibitem{Fan:Li:01}
Fan, J. and  Li, R.,(2001). 
\newblock  Variable selection via nonconcave penalized likelihood and its oracle properties.
\newblock  {\it Journal of the American Statistical Association}, \textbf{96}(456),  1348-1360.
\bibitem{Koenker:05}
Koenker, R.,(2005). 
\newblock \textit{Quantile Regression}, 
\newblock{Cambridge University Press}.
\bibitem{Xu:Ying:10}
Xu, J. and Ying, Z.,(2010). 
\newblock Simultaneous estimation and variable selection in median regression using Lasso-type penalty.
\newblock {\it Annals of the Institute of Statistical Mathematics}, \textbf{62},  487-514.
\bibitem{Wang:12}
Wang L.,(2013). 
\newblock $L_1$ penalized LAD estimator for high dimensional linear regression.
\newblock {\it Journal of Multivariate Analysis},  DOI 10.1016/j.jmva.2013.04.001.
\bibitem{Wu:Liu:09}
Wu, Y. and Liu, Y.,(2009). 
\newblock Variable selection in quantile regression.
\newblock {\it Statistica Sinica}, \textbf{19},  801-817.
\bibitem{Zou:06}
Zou, H.,(2006). 
\newblock  The adaptive Lasso and its oracle properties.
\newblock  {\it Journal of the American Statistical Association}, \textbf{101}(476),  1418-1428.
\end{thebibliography}
\end{document}